\documentclass[12pt,a4paper]{article}
\usepackage{amsthm}
\usepackage{amsfonts}
\usepackage{amsmath}
\usepackage{amssymb}
\usepackage{graphicx}
\usepackage{float}
\usepackage{caption}
\usepackage{subcaption}
\usepackage{mathrsfs}
\usepackage{hyperref}
\usepackage{graphicx}
\usepackage{tikz-cd}
\usetikzlibrary{babel}
\usepackage{comment}
\usepackage{geometry}
\usepackage{relsize}
\usepackage{mathtools}
\usepackage{adjustbox}
\usepackage{xcolor}
\newtheorem{nota}{Remark}
\usepackage{bm}

\title{Reduced-Order Time Splitting for Navier-Stokes with Open Boundaries}

\author{
 M. Aza\"iez\thanks{Bordeaux University, Bordeaux INP and I2M (UMR CNRS 5295), France. {\tt azaiez@u-bordeaux.fr}}, \quad
T. Chac\'on Rebollo \thanks{Departamento EDAN \& IMUS, Universidad de Sevilla, Spain.  {\tt chacon@us.es}}, \quad
C. N\'u\~nez Fern\'andez   \thanks{Departamento EDAN \& IMUS, Universidad de Sevilla, Spain. {\tt cnfernandez@us.es}}, \quad
S. Rubino \thanks{Departamento EDAN \& IMUS, Universidad de Sevilla, Spain.  {\tt samuele@us.es}}}

\date{\today}

\begin{document}
\maketitle

\begin{abstract}
In this work, we propose a Proper Orthogonal Decomposition-Reduced Order Model (POD-ROM) applied to time-splitting schemes for solving the Navier-Stokes equations with open boundary conditions. In this method, we combine three strategies to reduce the computing time to solve NSE: time splitting, reduction of the computational domain through non-standard treatment of open boundary conditions and reduced order modelling. To make the work self-contained, we first present the formulation of the time-splitting scheme applied to the Navier-Stokes equations with open boundary conditions, employing a first-order Euler time discretization and deriving the non-standard boundary condition for pressure. Then, we construct a Galerkin projection-based ROM using POD with two different treatments of the pressure boundary condition on the outlet. We propose a comparative performance analysis between the standard projection-based POD-ROM (fully intrusive) and a hybrid POD-ROM that combines a projection-based approach (intrusive) with a data-driven technique (non-intrusive) using Radial Basis Functions (RBF). We illustrate this comparison through two different numerical tests: the flow in a bifurcated tube and the benchmark numerical test of the flow past cylinder, numerically investigating the efficiency and accuracy of both ROMs.

\end{abstract}

\section{Introduction}

Extensive research has been conducted to reduce the computational cost in the simulation of incompressible flows. Many approaches focus on decoupling velocity from pressure using projection techniques, such as those proposed in \cite{chorin,chorin68,Temam69}. These techniques fall under the class of time-splitting methods. Many applications require solving the Navier-Stokes equations under open boundary conditions to avoid infinite domains. To address this fact, \cite{MejdiTimeSplit} introduces improvements for open and traction boundary conditions in time-splitting methods. In particular, the authors propose a method where the pressure boundary condition is derived by solving a low-dimension problem at the outflow boundary, allowing of reduction the computational domain while ensuring proper flow development.

On the other hand, Reduced-Order Models (ROMs) applied to the incompressible Navier-Stokes equations have been widely in the recent years studied due to their broad industrial applications \cite{Stabile,giovanni-gianluigi,REBOLLO2000301,iliescu2014variational}. However, when applying ROMs to fully coupled velocity-pressure formulations, a fundamental challenge arises: how to accurately recover or directly obtain the pressure approximation, given the weakly divergence-free nature of the reduced velocity space. Several strategies have been proposed to tackle this issue, including Local Projection Stabilization (LPS) methods \cite{samu,samujulia} and supremizers \cite{Ballarin2015supre}, among others \cite{Azaiez24,rebollo2022error,IliescuPressure}. Despite these advances, relatively few studies have been explored ROMs for problems involving pressure boundary conditions, mainly due to the difficulty of preserving these conditions at the reduced-order level \cite{RozzaHybrid}.

ROM methods can be broadly classified into projection-based (intrusive) approaches \cite{IliescuWang2014,janes2025new,rebollo2022error} and data-driven (non-intrusive) approaches \cite{giovanni-gianluigi,RozzaHybrid,Nunez2024,BirgulTraian}. Projec\-tion-based ROMs derive a reduced-order formulation directly from the governing equations. Typically, Proper Orthogonal Decomposition (POD) combined with Galerkin projection is employed to map the high-dimensional system onto a lower-dimensional space. These methods offer strong physical interpretability but require direct access to and modification of the full-order model (FOM), which can be restrictive in certain applications. In contrast, data-driven ROMs do not rely on direct access to the governing equations. Instead, they infer the reduced dynamics from data, making them particularly useful for black-box models or experimental datasets. Common techniques include Radial Basis Functions (RBF) and Artificial Neural Networks (ANN). While non-intrusive ROMs offer greater flexibility, their accuracy and generalization capabilities strongly depend on the quality and representativeness of the training data.

Although ROMs solution of incompressible fluid flows have been extensively studied, their application to time-splitting methods remains limited. Notable contributions
include \cite{Li22}, where Li et al. propose a POD-ROM based on the classical Chorin-Temam projection, and \cite{Azaiez2025}, which presents a POD-ROM for the Goda time-splitting scheme.
 
In the present work, we propose a POD-ROM applied to time-splitting methods for the Navier-Stokes problem with open boundary conditions. In this way, we exploit the combination of three strategies to reduce the computing time to solve the NSE: time splitting, reduction of computational domain through non-standard treatment of open boundary conditions and reduced order modelling. After explicitly detailing the fully intrusive POD-ROM, we introduce a hybrid approach that combines a data-driven approximation using Radial Basis Functions (RBF) to approximate the pressure boundary condition with an intrusive method for velocity computation. We investigate the efficiency and accuracy of both approaches in two different numerical tests: flow in a bifurcated tube and flow around a cylinder.
 
Within this framework, the main contribution of the paper is to devise a POD-ROM to solve NSE in pipe-like domains, that reaches computational speed-ups of over $1500$ compared to the FOM, with error levels of nearly $1\%$.

The structure of the paper is as follows.  Section \ref{SectionFOM} presents the application of a time-splitting scheme to the first-order Euler time discretization of the variational formulation. This includes the derivation of the pressure boundary condition and the description of a Finite Element (FE) approximation, leading to the fully discrete Full-Order Model (FOM) used to generate snapshots for the Reduced-Order Model (ROM). Section \ref{ROMtimesplit} explicitly details the standard Galerkin projection-based POD-ROM model and the hybrid POD-ROM model, which integrates the non-intrusive data-driven RBF technique along with its adaptation for time extrapolation. Section \ref{NumericalStudies} describes numerical results for two different numerical tests: parametrized flow in a bifurcated tube and unsteady flow around a cylinder, comparing the performance of both ROM approaches analyzed. Finally, Section \ref{Conclusion} summarizes the main conclusions and outlines potential future work.

\section{Time splitting approximation}\label{SectionFOM}
In this section, we introduce the proposed Full Order Model (FOM) for the incompressible Navier-Stokes Equations (NSE), and describe both its time discretization and the spatial approximation.

We consider the incompressible evolution NSE problem:
\begin{equation}\label{NavierStokesProblem}
\left\{ \begin{array}{rll}
    \bm{u}_t + (\bm{u} \cdot \nabla) \bm{u} - \nu \Delta \bm{u}  + \nabla p  & = \bm{f} & \text{ in } \Omega \times (0,T], \\
    \nabla \cdot \bm{u}  & = 0 & \text{ in } \Omega \times (0,T], \\
    \bm{u} & = \bm{0} & \text{ on } \Gamma_D\times (0,T], \\
    (\nu \nabla \bm{u} - p\bm{I}) \cdot\bm{n} & = \bm{0} & \text{ on } \Gamma_N\times (0,T], \\
    \bm{u}(x,0) & = \bm{u}_0(x) & \text{ in } \Omega,
\end{array} \right.
\end{equation}
where $\Omega \subset \mathbb{R}^{d}, \, d \in \{2,3\},$ is a bounded polyhedral domain with a Lipschitz-continuous boundary $ \Gamma = \partial \Omega = \overline{\Gamma}_D \cup \overline{\Gamma}_N,$ $[0,T]$ the time interval, and $\bm{u}_0(x)$ be the initial condition. In the above equations \eqref{NavierStokesProblem}, $\bm{u}$ is the velocity field, $p$ is the pressure of the incompressible fluid, $\nu > 0$ is the kinematic viscosity, and $\bm{f}$ is the forcing term. For the sake of simplicity, we impose homogeneous Dirichlet boundary conditions on $\Gamma_D $ and do-nothing (open) boundary conditions on $\Gamma_N.$ More general inflow boundary conditions may be taken into account by standard lifting techniques for NSE.

\subsection{Time splitting scheme for open boundary conditions}\label{SectionTimeSplit}

In this section, we describe the time discretization of problem \eqref{NavierStokesProblem} by time splitting scheme that includes a a non-standard pressure-correction scheme for open boundary conditions, already introduced in \cite{MejdiTimeSplit} for a finite volume FOM. We assume that $\Omega$ is a $2D$ channel-like flow, with inflow boundary $\Gamma_D$ and outflow boundary $\Gamma_N.$ This boundary is assumed to be a vertical straight segment of height $2L$ from $x_2 = -L$ to $x_2 = L$. Extension to a more general definition of $\Gamma_N$ can be found in \cite{MejdiTimeSplit}, Section 4.1. This is a  semi-discrete (in time) approximation of \eqref{NavierStokesProblem} that we describe formally, without explicit mention of the regularity needed by the different unknowns that appear in it.\\

Denoting $\bm{u}^n \approx \bm{u}(\cdot,t^n)$ and $p^n \approx p(\cdot,t^n)$ approximations of the solution of \eqref{NavierStokesProblem} at a certain time $t^n = n \Delta t, \,  n \in \{1, \ldots, N\},$ with $\Delta t = {T}/{N},$ The time splitting scheme that we consider in this work operates in two stages. We obtain $\bm{u}^{n+1}$ and $p^{n+1}$ by solving: 

\noindent \textbf{Step 1.} Convection-diffusion problem:
\emph{find $\bm{\tilde{u}}^{n+1}$ such that} 
\begin{equation}\label{FinalStep1}
    \begin{array}{rll}
        \dfrac{\bm{\tilde{u}}^{n+1}-\bm{u}^n}{\Delta t} + (\bm{u}^{n} \cdot \nabla) \bm{\tilde{u}}^{n+1} - \nu \Delta \bm{\tilde{u}}^{n+1} + \nabla p^n & = \bm{f}^{n+1} & \text{ in } \Omega,  \\
         \bm{\tilde{u}}^{n+1} & = \bm{0} & \text{ on }  \Gamma_D, \\
         (\nu \nabla \bm{\tilde{u}}^{n+1} - p^{n}\bm{I})\cdot\bm{n} & = \bm{0}& \text{ on } \Gamma_N.
\end{array}
\end{equation}

\noindent \textbf{Step 2.} Open boundary pressure BC problem: 
\emph{find $\hat{\phi}^{n+1}$ defined on $\Gamma_N$ such that}
\begin{equation}\label{FinalStep1.5}
\left \{\begin{array}{rll}
    \left( \Delta t \partial_{x_2}^2 - \dfrac{1}{\nu}\right)  \hat{\phi}^{n+1} & = \nabla \cdot \bm{\tilde{u}}^{n+1} & \text{ on } \Gamma_N, \\
    \partial_{\bm{n}} \hat{\phi}^{n+1}(\pm L) & = 0.&  \\
\end{array} \right. 
\end{equation}

\noindent \textbf{Step 3.} Pressure-continuity correction problem: 
\emph{find $\phi^{n+1}$ such that }
\begin{equation}\label{FinalStep2}
    \begin{array}{rll}
        \Delta t \Delta \phi^{n+1}  & = \nabla \cdot \bm{\tilde{u}}^{n+1} & \text{ in } \Omega,  \\
        \partial_{\bm{n}} \phi^{n+1} & = 0 & \text{ on } \Gamma_D, \\
        \phi^{n+1} & = \hat{\phi}^{n+1} & \text{ on } \Gamma_N.
    \end{array}
\end{equation}

\medskip

\noindent \textbf{Step 4.} Updating velocity and pressure problem: 
\emph{find $\bm{u}^{n+1}$ and $p^{n+1}$ such that}
\begin{equation}\label{FinalStep3pres}
    p^{n+1} = p^{n} + \phi^{n+1} \text{ in } \Omega,
\end{equation}
\begin{equation}\label{FinalStep3vel}
\bm{u}^{n+1} = \bm{\tilde{u}}^{n+1} - \Delta t \nabla \phi^{n+1} \text{ in } \Omega.
\end{equation}

\subsubsection*{Pressure outflow boundary conditions}
The boundary condition for
$\phi^{n+1}$ on $\Gamma_N$ is deduced in \cite{MejdiTimeSplit} by derivation of the Helmholtz-Hodge decomposition of
$\tilde{u}^{n+1}$ in \eqref{FinalStep2}, as also used \cite{jin1993nonreflecting,kirkpatrick2008open}  for other kinds of outlet boundary conditions
(non-reflecting or Neumann). As mentioned in \cite{guermond}, the natural choice to impose $\phi^{n+1}$ = 0 on $\Gamma_N$ would lead to numerical issues, as it constrains the pressure increment and can
prevent its proper evolution near the boundary. The extension to boundaries $\Gamma_N$ with non-zero curvatures can be found in Section $4.1$ of \cite{MejdiTimeSplit}.

\subsection{FE space approximation}\label{FETimeSplit}
In this section, we state the variational formulation of the considered time-splitting method and its Finite Element (FE) discretization. This will be the FOM to solve NSE \eqref{NavierStokesProblem} that we consider in this paper.

Let us consider the following spaces: 
\begin{itemize}
\item $ \bm{H}^1_D(\Omega) = \{ \bm{\tilde{v}} \in \bm{H}^1(\Omega) : \bm{\tilde{v}} = \bm{0} \text{ on } \Gamma_D \}.$
   \item $ \bm{H}(div,\Omega) =  \{ \bm{v} \in \bm{L}^2(\Omega) : \, \nabla \cdot \bm{v} \in L^2(\Omega),\bm{v}\cdot\bm{n}=0 \text{ on } \Gamma_D \}. \, $
   \item $H^1_N(\Omega) = \{ q \in H^1(\Omega) : q = 0 \text{ on } \Gamma_N \}.$
\end{itemize}

The variational formulation of the FOM scheme described in Section \ref{SectionTimeSplit} reads:

\noindent \textbf{Step 1.} Prediction-diffusion problem:  \emph{given $\bm{u}^n \in \bm{H}(div,\Omega), p^n \in L^2(\Omega)$ and $\bm{f}^{n+1} \in \bm{H}^{-1},$ find $\bm{\tilde{u}}^{n+1} \in \bm{H}^1_D(\Omega)$ such that, for all $\bm{\tilde{v}} \in \bm{H}^1_D(\Omega) $}
\begin{equation}\label{Step1FV}
        \dfrac{1}{\Delta t}(\bm{\tilde{u}}^{n+1}-\bm{u}^n,\bm{\tilde{v}}) + (\bm{u}^{n}\cdot\nabla\bm{\tilde{u}}^{n+1},\bm{v}) + \nu (\nabla \bm{\tilde{u}}^{n+1}, \bm{\tilde{v}}) - (p^n,\nabla \cdot \bm{\tilde{u}}^{n+1}) = \langle \bm{f}^{n+1},\bm{\tilde{v}} \rangle.
\end{equation}

\noindent \textbf{Step 2.} Open Boundary pressure BC problem:  
\emph{find $\hat{\phi}^{n+1} \in H^1(\Gamma_N)$ such that} 
\begin{equation}\label{Step15FV}
\nu \Delta t(\partial_{x_2} \hat{\phi}^{n+1},\partial_{x_2}\hat{q})_{\Gamma_N} + (\hat{\phi}^{n+1},\hat{q})_{\Gamma_N} = -\nu (\nabla \cdot \bm{\tilde{u}}^{n+1},\hat{q})_{\Gamma_N} \hspace{0.5cm} \forall \hat{q} \in H^1(\Gamma_N).
\end{equation}

\noindent  \textbf{Step 3.} Pressure-continuity correction problem: \emph{find $\phi^{n+1} \in H^1_N(\Omega)$ such that}
\begin{equation}\label{Step2FV}
    (\nabla \phi^{n+1}, \nabla \tilde{q}) = -\dfrac{1}{\Delta t} (\nabla \cdot \bm{\tilde{u}}^{n+1},\tilde{q}) - (\nabla \tilde{\phi}^{n+1}, \nabla \tilde{q}) \hspace{0.5cm} \forall \tilde{q} \in H^1_N(\Omega),
\end{equation}
where $\tilde{\phi}^{n+1} \in H^1(\Omega)$ is an extension of the solution of \eqref{Step15FV} to $\Omega,$ then  $\phi^{n+1}$ verifies homogeneous Dirichlet boundary condition in $\Gamma_{N}.$\\

\noindent \textbf{Step 4. } Updating velocity and pressure problem:  
\emph{find $p^{n+1} \in L^2(\Omega)$ such that}
\begin{equation}\label{Step3PresFV}
    p^{n+1} = p^n + \overline{\phi}^{n+1},
\end{equation}
\emph{and $\bm{u}^{n+1} \in  \bm{H}(div,\Omega)$ such that}
\begin{equation}\label{Step3VelFV}
    (\bm{u}^{n+1},\bm{v}) = (\bm{\tilde{u}}^{n+1},\bm{v}) - \Delta t(\nabla \overline{\phi}^{n+1},\bm{v}), \hspace{0.5cm} \forall \bm{v} \in \bm{H}(div,\Omega), 
\end{equation}
where $\overline{\phi}^{n+1} = \phi^{n+1} + \tilde{\phi}^{n+1}.$

To derive a FE approximation of the variational scheme \eqref{Step1FV}-\eqref{Step3VelFV}, we apply a Galerkin projection onto a suitable finite-dimensional subspaces of the corresponding functional spaces,
\begin{equation}\label{FES}
    \bm{X}_h \subset \bm{H}(div,\Omega), \hspace{0.25cm} \bm{\tilde{X}}_h \subset \bm{H}^1_D(\Omega),  \hspace{0.25cm} \tilde{Q}_h \subset H^1_N(\Omega), \hspace{0.25cm} \hat{Q}_h \subset H^1(\Gamma_N),\hspace{0.25cm}  Q_h \subset L^2(\Omega).
\end{equation}

The fully discrete scheme is then obtained by seeking, at each time step, the approximate solutions $\bm{u}_h^{n+1}, \bm{\tilde{u}}_h^{n+1}, p_h^{n+1},\phi^{n+1}_h$ and $\hat{\phi}^{n+1}_h$ in the above FE spaces that satisfy the weak formulation \eqref{Step1FV}--\eqref{Step3VelFV} in a Galerkin sense. For brevity we do not make explicit this discretization as it is straighforward. It follows by changing every space in \eqref{Step1FV}-\eqref{Step3VelFV} by its discrete counterpart, following \eqref{FES}.

\section{POD-ROM method}\label{ROMtimesplit}

In this section, we present the POD-ROM approximation to the fully discrete time-splitting method  introduced in Section \ref{SectionFOM}.

We describe the offline phase for the standard POD-ROM projection-based method and then, the online phase. 

\subsection{Offline Phase}\label{OfflinePhase}
We start from the solution of the FOM problem, that provides snapshots for velocity $\{ \bm{u}_h^{t_i} \}_{i=1}^N,$ predicted velocity $\{ \bm{\tilde{u}}_h^{t_i} \}_{i=1}^N,$ pressure $\{p_h^{t_i} \}_{i=1}^N,$ pressure-continuity correction $\{ \phi_{h}^{t_i}\}_{i=1}^N$ and outflow pressure BC $\{ \hat{\phi}_{h}^{t_i}\}_{i=1}^N.$

The POD method constructs low-dimensional spaces to approximate each variable of the scheme, whose bases we respectively denoted as $\{\bm{\varphi}_1, \ldots, \bm{\varphi}_{r_{\bm{u}}}  \},$ $ \{ \bm{\tilde{\varphi}}_1, \ldots, \bm{\tilde{\varphi}}_{r_{\bm{\tilde{u}}}} \},$ $\{ \psi_1, \ldots, \psi_{r_p} \}, $ $\{ \gamma_1, \ldots, \gamma_{r_{\phi}}\} $ and $\{ \hat{\gamma}_1, \ldots, \hat{\gamma}_{r_{\hat{\phi}}}\}$. These spaces are designed to optimally approximate the given snapshots sets in a least-squares sense with respect to Hilbertian norms \cite{kunisch2001galerkin}. In what follows, we consider POD spaces computed using $L^2$ scalar products for velocity and predicted velocity fields, while pressure spaces are constructed using $H^1$   scalar products to account for their inherent continuity and smoothness requirements. We assume that the selected POD basis retains a high percentage of the total energy carried by the initial snapshots. We thus consider the following  reduced spaces for the POD setting:

\begin{itemize}
    \item $\bm{\mathcal{U}}_{r_{\bm{u}}} = \text{ span } \{\bm{\varphi}_1, \ldots, \bm{\varphi}_{r_{\bm{u}}}  \}$ \,\,(reduced velocity  space)
    \item $\bm{\tilde{\mathcal{U}}}_{r_{\bm{\tilde{u}}}} = \text{ span } \{ \bm{\tilde{\varphi}}_1, \ldots, \bm{\tilde{\varphi}}_{r_{\bm{\tilde{u}}}} \}$ \,\,(reduced predicted velocity  space)
    \item $\mathcal{Q}_{r_p} = \text{ span } \{ \psi_1, \ldots, \psi_{r_p} \},$ \,\,(reduced pressure  space)
    \item $\tilde{\mathcal{Q}}_{r_\phi} = \text{ span } \{ \gamma_1, \ldots, \gamma_{r_{\phi}} \}$ \,\,(reduced pressure-continuity correction  space)
    \item  $\hat{\mathcal{Q}}_{r_{\hat{\phi}}} = \text{ span } \{ \hat{\gamma}_1, \ldots, \hat{\gamma}_{r_{\hat{\phi}}} \}$\,\, (reduced pressure-continuity correction boundary  space)
\end{itemize}

\subsection{Online Phase }

In this section, we propose a ROM approximation of the FE approximation of the splitting method \eqref{Step1FV}-\eqref{Step3PresFV} using Galerkin projection. We approximate the FE solution $(\bm{\tilde{u}}_h^{n+1}, \hat{\phi}_h^{n+1}, \phi_h^{n+1}, \bm{u}_h^{n+1}, p_h^{n+1})$  by its reduced-order counterpart $(\bm{\tilde{u}}_r^{n+1}, \hat{\phi}_r^{n+1}, \phi_r^{n+1}, $ $ \bm{u}_r^{n+1}, p_r^{n+1})$ by 
\begin{equation}
    \bm{u}(x,t) \approx \bm{u}_{r_{\bm{u}}}(x,t) = \sum_{i=1}^{r_{\bm{u}}} a_i(t) \bm{\varphi}_i(x) \in \bm{\mathcal{U}}_{r_{\bm{u}}},
\end{equation}

\begin{equation} 
    \bm{\tilde{u}}(x,t) \approx \bm{\tilde{u}}_{r_{\bm{\tilde{u}}}} (x,t) 
    = \sum_{i=1}^{r_{\bm{\tilde{u}}}} \tilde{a}_i(t)\bm{\tilde{\varphi}}_i(x) \in \bm{\tilde{\mathcal{U}}}_{r_{\bm{\tilde{u}}}},
\end{equation}

\begin{equation}
    p(x,t) \approx p_{r_p}(x,t) = \sum_{i=1}^{r_p} b_i(t)\psi_i(x) \in \mathcal{Q}_{r_{p}},
\end{equation}

\begin{equation}
    \phi(x,t) \approx \phi_{r_\phi}(x,t) = \sum_{i=1}^{r_\phi} c_i(t)\gamma_i(x) \in \tilde{\mathcal{Q}}_{r_{\phi}},
\end{equation}

\begin{equation}
    \hat{\phi}(x,t) \approx \hat{\phi}_{r_{\hat{\phi}}}(x,t) = \sum_{i=1}^{r_{\hat{\phi}}} \hat{c}_i(t)\hat{\gamma}_i(x) \in \hat{\mathcal{Q}}_{r_{\hat{\phi}}},
\end{equation}

\noindent where $\{a_i(t)\}_{i=1}^{r_u},\{\tilde{a}_i(t)\}_{i=1}^{r_{\tilde{u}}}, \{b_i(t)\}_{i=1}^{r_p} $, $\{c_i(t)\}_{i=1}^{r_{\phi}}$ and $\{\hat{c}_i(t)\}_{i=1}^{r_{\hat{\phi}}} $are the sought time-varying coefficients for each field, respectively. Note that the dimension of the corresponding vectors is much smaller than the number of degrees of freedom in a full order simulation.

Projecting the equations of the FOM splitting method \eqref{Step1FV}-\eqref{Step3VelFV} onto their respective POD spaces and assuming that we know an approximation of the velocity $\bm{u}_h$ and the pressure $p_h$ at time $t=t_1,$ the fully discrete space-time formulation of the Galerkin POD-ROM is derived as follows:\\

\noindent \textbf{Step 0.} Initialization. 
\begin{itemize}
    \item  $\bm{u}^1_r = \sum_{i=1}^{r_{\bm{u}}} (u^{t_1}_h,\bm{\varphi_i})\bm{\varphi}_i,$
    \item $p^1_r = \sum_{i=1}^{r_p} (p^{t_1}_h,\psi_i)\psi_i.$
\end{itemize}

For $n=1,\ldots,N:$ \\

\noindent \textbf{Step 1.} Prediction-diffusion reduced order problem:  
given $\bm{u}_n^r \in \bm{\mathcal{U}}_{r_u}$ and $p_r^n \in \mathcal{Q}_{r_p}$, find $\bm{\tilde{u}}_r^{n+1} \in \bm{\tilde{\mathcal{U}}}_{r_{\tilde{u}}}$ such that, for all $  \bm{\tilde{\varphi}} \in \bm{\tilde{\mathcal{U}}}_{r_{\bm{u}}}$
\begin{equation}\label{ROMStep1}
    \begin{array}{rl}
        (\dfrac{\bm{\tilde{u}}_r^{n+1}-\bm{u}_r^n}{\Delta t},\bm{\tilde{\varphi}}) + (\bm{u}_r^n \cdot \nabla \bm{\tilde{u}}^{n+1}_r,\bm{\tilde{\varphi}})  + \nu (\nabla \bm{\tilde{u}}_r^{n+1}, \nabla \bm{\tilde{\varphi}}) - (p_r^n, \nabla \cdot \bm{\tilde{\varphi}}) & = \langle \bm{f}^{n+1},\bm{\tilde{\varphi}} \rangle.  \\
    \end{array}
\end{equation}

\noindent \textbf{Step 2.} Open Boundary pressure BC reduced order problem:  
given $\bm{\tilde{u}}_r^{n+1} \in \bm{\tilde{\mathcal{U}}}_{r_{\tilde{u}}}, $ find $\hat{\phi}_r^{n+1} \in \hat{\mathcal{Q}}_{r_{\hat{\phi}}}$ such that
\begin{equation}\label{ROMStep15}
\nu \Delta t(\partial_{x_2} \hat{\phi}^{n+1}_r,\partial_{x_2}\hat{\gamma})_{\Gamma_N} + (\hat{\phi}^{n+1}_r,\hat{\gamma})_{\Gamma_N} = -\nu (\nabla \cdot \bm{\tilde{u}}_r,\hat{\gamma})_{\Gamma_N} \hspace{0.5cm} \forall \hat{\gamma} \in \hat{\mathcal{Q}}_{r_{\hat{\phi}}}.
\end{equation}

\noindent \textbf{Step 3.} Pressure-continuity correction reduced order problem: 
given $\bm{\tilde{u}}_r^{n+1} \in \bm{\tilde{\mathcal{U}}}_{r_{\tilde{u}}}$ and $\hat{\phi}_r^{n+1} \in \hat{\mathcal{Q}}_{r_{\hat{\phi}}},$ find $\phi_r^{n+1} \in \tilde{\mathcal{Q}}_{r_{\phi}}$ such that
\begin{equation}\label{ROMStep2}
    (\nabla \phi^{n+1}_r, \nabla \gamma) = -\dfrac{1}{\Delta t} (\nabla \cdot \bm{\tilde{u}}^{n+1}_r,\gamma) - (\nabla \tilde{\phi}^{n+1}_r, \nabla \tilde{q}_h) \hspace{0.5cm} \forall \gamma \in \tilde{\mathcal{Q}}_{r_{\gamma}},
\end{equation}
where $\tilde{\phi}^{n+1}_r$ is the extension (constant with respect $y$) of the solution of \eqref{ROMStep15} to the full domain.\\

\noindent \textbf{Step 4.} Updating pressure and velocity   reduced order problem: 
 given $p^n_r \in \mathcal{Q}_{r_p},$ $\phi^{n+1}_r \in \tilde{\mathcal{Q}}_{r_{\phi}}$ and $ \tilde{\phi}_r^{n+1} \in \hat{\mathcal{Q}}_{r_{\hat{\phi}}},$ find $p_r^{n+1} \in \mathcal{Q}_{r_p}$ such that 
\begin{equation}\label{ROMStep3Pres}
    \begin{array}{rl}
         (p^{n+1}_r, \psi) & = (p^n_r,\psi) + (\overline{\phi}^{n+1}_r,\psi)\quad  \forall \psi \in \mathcal{Q}_{r_p}.  
    \end{array}
\end{equation}
And given $\bm{\tilde{u}}^n_r \in \bm{\tilde{\mathcal{U}}}_{r_{\tilde{u}}}$, $\phi^{n+1}_r \in \tilde{\mathcal{Q}}_{r_{\phi}} $ and $\tilde{\phi}_r^{n+1} \in \hat{\mathcal{Q}}_{r_{\hat{\phi}}}$, find $\bm{u}_r^{n+1} \in \bm{\mathcal{U}}_{r_{\bm{u}}}$ such that:
\begin{equation}\label{ROMStep3Vel}
    \begin{array}{rl}
         (\bm{u}^{n+1}_r,\bm{\varphi}) & = (\bm{\tilde{u}}^{n+1}_r,\bm{\varphi}) - \Delta t (\nabla \overline{\phi}^{n+1}_r, \bm{\varphi}) \quad\forall \bm{\varphi} \in \bm{\mathcal{U}}_{r_{\bm{u}}}, 
    \end{array}
\end{equation}
where $\overline{\phi}^{n+1}_r = \phi^{n+1}_r + \tilde{\phi}_r^{n+1}.$

\medskip 

Note that all the involved matrices corresponding to the reduced problems are   precomputed during the offline stage. 

\section{Hybrid POD-ROM method}\label{HybridROM}
In this section, we present a hybrid POD-ROM strategy that combines a standard intrusive Galerkin projection approach with a non-intrusive/data-driven technique based on RBF method. The goal is to enhance the standard POD-ROM framework by incorporating a pressure-continuity correction computed through machine learning techniques.

We begin by introducing the standard RBF interpolation method, applied to multi-parameter settings including time and a selected physical parameter. This method is particularly effective in interpolation regimes but lacks robustness when extrapolating in time. To address this limitation, we then describe an RBF-based extrapolation strategy specifically designed for time periodic. Finally, we outline the online phase of the hybrid POD-ROM, where the intrusive ROM is coupled with the non-intrusive correction.

\subsection{Radial Basis Method for the computation of pressure-continuity correction}\label{RBFsection}

We build in this section an RBF approximation of unknowns $\phi^{n+1}_r$ and $\hat{\phi}^{n+1}_r,$ respectively solution of Steps \eqref{ROMStep2} and \eqref{ROMStep15}. 
The method will be explained only for $\phi_r^{n+1},$ its adaptation to $\hat{\phi}^{n+1}_r$ is straightforward. We have considered a multiparametric dependence on time and physical parameter (Reynolds number).\\

In this framework, the RBF strategy starts from the POD approximation of $\phi(x,t; \mu)$ in the form:
\begin{equation}\label{RBFphi}
    \phi(x,t;\mu) \approx \phi_{r_{\phi}} = \sum_{i=1}^{r_{\phi}} c_i(t;\mu)\gamma_i(x),
\end{equation}
where $ \{\gamma_i \}_{i=1}^{r_{\phi}}$ is the POD basis associated to the pressure-continuity correction and $\{c_i(t;\mu)\}_{i=1}^{r_{\phi}}$ is the parameter-depending coefficients vector. Our strategy is to approximate each of these coefficients by the RBF, adapting the procedure introduced in \cite{giovanni-gianluigi}. 

Let us introduce the following notation: 
\begin{itemize}
\item $X_{\mu,t} = \{ \mu_1,\ldots,\mu_M\} \times \{t_1,\ldots,t_N\}$ is the cartesian product of the discretized parameter set and the set of time instances in which pressure correction snapshots were taken.  $x^{i}_{\mu,t}$ is the i-th member of $X_{\mu,t}.$
\item $z^* = (t^*,\mu^*),$ where  $t^{*}$ is the time instant at which the ROM solution is sought and $\mu^*$ is the new parameter considered in the online stage of the ROM (within the range $[\mu_1,\mu_M],$ but different from the parameters used in the offline stage).  
\item $S = \{\phi(x,t_1;\mu_1),\ldots,\phi(x,t_N;\mu_M)\} \in \mathbb{R}^{N^h_{\phi} \times \tilde{M}}$ is the set of snapshots, where $\tilde{M} = N \cdot M$ and $N_h^\phi$ is the number of degrees of freedom of the pressure-continuity correction. 
\item $\tilde{c}_{r,l} = (S^r,\gamma_l)$ is the projection of the $r$-th column of $S$ on the $l$-th POD mode. 
\end{itemize}

The interpolation will be carried out independently for each mode $\gamma_L$, so we can fix one pressure-continuity correction mode and focus on approximating its associated coefficient function $c_L(t;\mu).$ The procedure is based on the classical radial basis function interpolation:
\begin{equation}\label{RBFormula}
       c_L(z) = \sum_{j=1}^{\tilde{M}} w_{L,j} \xi (\|z-x_{\mu,t}^j\|_{L^2}) \text{ for } L = 1,2,\ldots, r_{\phi},
\end{equation}
where $z = (t,\mu), w_{L,j}$ are some appropriate weights and $\xi$ is the RBF function centered in $x_{\mu,t}^j.$ In contrast to the classical Gaussian RBF, we instead employ the regularized thin‑plate spline
\begin{equation*}
    \xi(r) = r^2 log(r+1),
\end{equation*}
which is radial and $C^{\infty},$ avoids the singularity at $r = 0$ of the standar $r^2log(r),$ and 
admits better numerical stability for scattered data in high‑dimensional parameter spaces (see Chapter 2 in \cite{ThesisManzoni}).
\medskip 

\subsubsection{Offline RBF Phase.} 
\begin{enumerate}
    \item We impose $c_L(x_{\mu,t}^i) = \sum_{j=1}^{\tilde{M}} w_{L,j}, \xi (\|x_{\mu,t}^i-x_{\mu,t}^j\|_{L^2}) = \tilde{c}_{i,L} $ for $i = 1, \ldots, \tilde{M}.$
    \item This corresponds to the following linear system: 
        \begin{equation}\label{RBFlinear}
            A^{\xi} w_L = Y_L,
        \end{equation}
        where 
        \begin{align}
            (A^{\xi})_{ij} = \xi (\|x^i_{\mu,t} - x^j_{\mu,t}\|_{L^2}), \\
            Y_L = \{\tilde{c}_{r,L} \}_{r=1}^{\tilde{M}}.
        \end{align}
        The linear system \eqref{RBFlinear} admits a unique solution under standard RBF assumptions as the centers $x_{\mu,t}^j$ are all different, and the kernel $\xi$ is positive defined (\cite{RBFSystem}). The weights $w_L$ can thus be computed in the offline phase, and then stored to be used in the online phase. 
\end{enumerate}

\subsubsection{Online RBF Phase.} 
Given a new parameter-time pair $z^* = (t^*,\mu^*),$ the corresponding coefficient is evaluated as:
\begin{equation}
    c_L(z^*) = \sum_{j=1}^{\tilde{M}} w_{L,j} \xi(\|z^* - x^j_{\mu,t}\|_{L^2}).
\end{equation}
This is repeated for each mode $L = 1, \ldots, r_{\phi}$ to recover all reduced coefficients.

\subsection{RBF extrapolation}\label{RBFextrapolation}

We also consider an alternative strategy of approximation of the pressures $\phi$ and $\tilde{\phi}$ adapted to time extrapolation. The idea comes from the fact in Step 2 of the time splitting method, detailed in Section \ref{SectionTimeSplit}, $\phi^{n+1}$ is computed through $\nabla \cdot \tilde{u}^{n+1}$, that is, via the predicted velocity $\tilde{u}^{n+1}$. Therefore, in the ROM, we propose to compute the coefficients $\bm{c}^{n+1}$ of $\phi^{n+1}$ through the coefficients of the predicted velocity, $\bm{\tilde{a}}^{n+1} = \{ {\tilde{a}}_i^{n+1} \}_{i=1}^{r_{\bm{\tilde{u}}}}$. This enables us to express $\bm{c}^{n+1}$ as:
\begin{equation}
    \bm{c}^{n+1} = \bm{c} ( t^{n+1};\mu) \approx \bm{c} (\bm{\tilde{a}}^{n+1}).
\end{equation}

To define the new RBF procedure, we consider $\bm{a}_j = \{a_{i,j}\}_{i=1,j=1}^{M,r_{\bm{\tilde{u}}}}$ instead of $X_{\mu,t},$ where each element is denoted by $a_{i,j} = (S^i_{\bm{\tilde{u}}},\bm{\tilde{\varphi}}_j).$ 

We now adapt the radial basis function interpolation formula \eqref{RBFormula} to this new setting:
\begin{equation}
    c_L(\bm{\tilde{a}}^{n+1}) = \sum_{j=1}^{\tilde{M}} w_{L,j} \xi (\|\bm{\tilde{a}}^{n+1} - \bm{a}_{j} \|_{L^2} ) \text{ for } L = 1,2, \ldots, r_{\phi}.
\end{equation}

\subsubsection{Offline RBF extrapolation Phase.} 
\begin{enumerate}
    \item We impose $c_L(\bm{a}_i) = \sum_{j=1}^{\tilde{M}} w_{L,j}, \xi (\|\bm{a}_i - \bm{a}_j\|_{L^2}) = \tilde{c}_{i,L} $ for $i = 1, \ldots, \tilde{M}.$
    \item This corresponds to the following linear system: 
        \begin{equation}
            A^{\xi} w_L = Y_L,
        \end{equation}
        where 
        \begin{equation}
            (A^{\xi})_{ij} = \xi (\|\bm{a}_i - \bm{a}_j\|_{L^2}).
        \end{equation}
\end{enumerate}

\subsubsection{Online RBF extrapolation Phase.} 
Given $\bm{\tilde{a}}^{n+1}$ the corresponding coefficient is evaluated as:
\begin{equation}
    c_L(\bm{\tilde{a}}^{n+1}) = \sum_{j=1}^{\tilde{M}} w_{L,j} \xi(\|\bm{\tilde{a}}^{n+1} - \bm{a}_j\|_{L^2}).
\end{equation}
This is repeated for each mode $L = 1, \ldots, r_{\phi}$ to recover all reduced coefficients.

\subsection{Online Phase for Hybrid POD-ROM method}\label{OnlinePhase}
Once assembled the offline matrices and detailed the RBF method for the computations of $\phi_r^{n+1}$ and $\hat{\phi}_r^{n+1}$ in Section \ref{RBFsection}, the online phase of the hybrid POD-ROM method is as follows.

\begin{itemize}
\item \textbf{Initialization.} The ROM initial conditions are obtained by POD spaces projection:
\begin{equation}
a_i^1 = (\bm{u}_{1h},\varphi_i) \hspace{0.5cm} i = 1,\ldots,r_{\bm{u}},
\end{equation}
\begin{equation}
b_i^1 = (\nabla p_{1h}, \nabla \psi_i) \hspace{0.5cm} i = 1,\ldots,r_{p}.
\end{equation}
Where $\bm{u}_{1h}$ and $p_{1h}$ are some approximations to $\bm{u}_1 = \bm{u}(\cdot,1)$ in $L^2-$norm belonging to $\bm{X}_h$ and and $p_1 = p(\cdot,1) $ in $H^1-$norm belonging to $Q_h,$ respectively.

\item \textbf{Iteration.} For $n = 1,\ldots,N$, perform one of the following ROM update schemes depending on the chosen approach:

\textbf{Step 1.} The same as the Galerkin Projection POD-ROM by solving the reduced linear system \eqref{ROMStep1}. \newline
\textbf{Step 2 and Step 3.} Reduced pressure-continuity correction: \newline 
Given $\bm{\tilde{a}}^{n+1}$ from Step 1, obtain $\bm{c}^{n+1}$ and $\bm{\hat{c}}^{n+1}$ by using the RBF described in Section \ref{RBFsection} (in the case of extrapolation by using the RBF described in Section \ref{RBFextrapolation}).\newline
\textbf{Step 4.} The same as the Galerkin Projection POD-ROM by solving the reduced linear systems \eqref{ROMStep3Pres}-\eqref{ROMStep3Vel}.
\end{itemize}

\section{Numerical studies}\label{NumericalStudies}
In this section, we present some numerical results to test and compare the POD-ROM and hybrid POD-ROM, introduced in Section \eqref{ROMtimesplit} and \eqref{HybridROM}, in terms of accuracy and efficiency. In the first numerical experiment, we employ the standard RBF interpolation method. In contrast, the second numerical experiment corresponds to the RBF extrapolation approach in time. The open-source FE software FreeFem \cite{Freefem} has been used to run the numerical experiments.

\subsection{Flow in a bifurcated tube}\label{Bifur}
\begin{figure}[ht!]
    \centering
    \includegraphics[width=\textwidth]{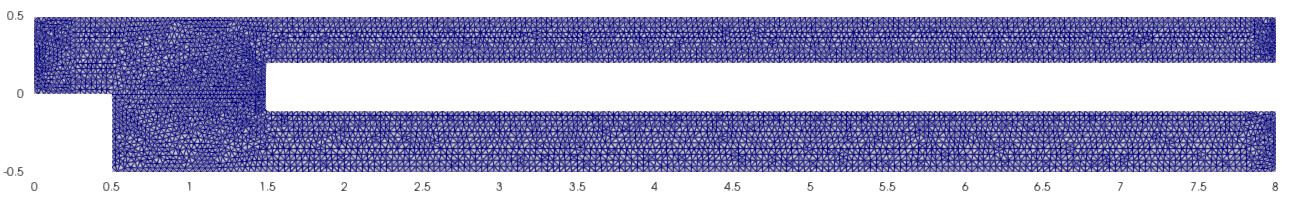}
    \caption{Example \ref{Bifur}: Computational grid of 2D flow in a bifurcated tube benchmark problem.}
    \label{fig:meshBifurcated}
\end{figure}

We consider the problem of 2D unsteady flow in a bifurcated tube as in \cite{MejdiTimeSplit,BifurcatedTube}. The computational domain is given by (see Figure \ref{fig:meshBifurcated}):
\begin{equation}
    \Omega = [0,8] \times [-0.5,0.5] \setminus \{[0,0.5] \times [-0.5,0] \cup [1.5,8] \times [-0.1,0.2] \}. 
\end{equation}

We impose at the inlet boundary $\{0\} \times \{0,0.5\}$ Dirichlet boundary conditions, by considering a Poiseuille flow with unitary influx. We also impose outflow boundary conditions at the outlet boundaries $\{8\} \times \{-0.5,-0.1\} \cup \{8\} \times \{0.2,0.5\}$ in order to test both ROMs detailed in Section \ref{ROMtimesplit}. No-slip boundary conditions are prescribed in the remaining walls.

To assess the computational results, we focus on the analysis of time evolution for the following three characteristic outputs of the flow: 

\begin{itemize}
    \item $E_{kin} = \frac12 \|\bm{u}\|_{L^2}^2 $(Kinetic Energy),
    \item $ Q = \displaystyle{\int}_R \bm{u} \cdot \bm{n} $ (Outflux),
    \item $CD = \frac12(\displaystyle{\int}_L (\bm{u} \cdot \bm{n})^2 - \displaystyle{\int}_{R}(\bm{u} \cdot \bm{n})^2) + (\displaystyle{\int}_L p - \displaystyle{\int}_R p)$ (Charge Drops),
\end{itemize}
where $L$ is the inflow boundary, and $R$ is the upper or lower outflow boundary.

For this test, we consider the time and the Reynolds number $Re = \dfrac{\bar{U} D}{\nu}$ as parameters, where $\bar{U} = 1\,m/s$ is the mean inflow velocity, $\nu$ the kinematic viscosity and $D = 1\,m$ is the total height of the bifurcated tube. The parameter range for the Reynolds number is $\mathcal{D} = [400,600]$ and for the time is $[5,7]s.$

\subsubsection{FOM and POD modes}
Our goal is to compare the fully-discrete time-splitting FOM described in Section \ref{SectionFOM} with both ROM counterparts described in Section \ref{ROMtimesplit} and \ref{HybridROM}. 

The FOM is the one described in Section \ref{SectionFOM} with a spatial discretisation using $\mathbb{P}^2$ FE for the velocity and predicted velocity and $\mathbb{P}^1$ FE for the pressure, pressure-continuity correction and pressure-continuity boundary correction on a relatively coarse computational grid (see Figure \ref{fig:meshBifurcated}), resulting in $5853$ and $22490$ degrees of freedom for pressure and velocity respectively. 

\begin{figure}[ht!]
\centering
\begin{subfigure}[b]{0.45\textwidth}
\includegraphics[width=\textwidth]{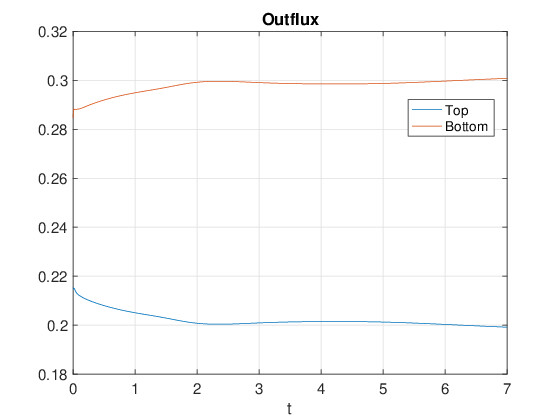}
\end{subfigure}
\begin{subfigure}[b]{0.45\textwidth}
\includegraphics[width=\textwidth]{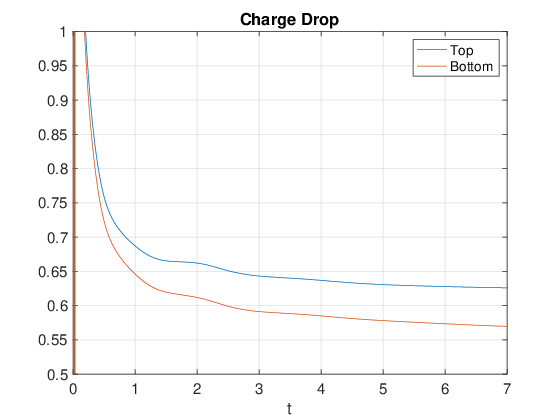}
\end{subfigure}
\begin{subfigure}[b]{0.45\textwidth}
\includegraphics[width=\textwidth]{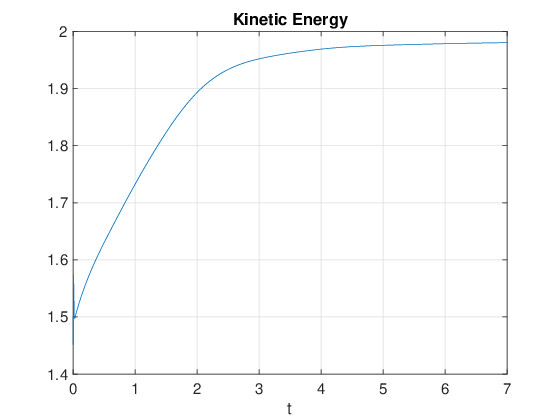}
\end{subfigure}
\caption{Example \ref{Bifur}: Temporal evolution of outflux (top left figure), charge drops (top right figure) and kinetic energy (bottom figure) for the FOM solution of $Re = 600.$}
\label{fig:QuantitiesBifurcated}
\end{figure}

The FOM simulation is initialized with zero velocity field and zero pressure and we have considered a time step $\Delta t = 10^{-2}.$ We have carried out the simulation up to a final time $T = 7s$, where the flow remains transient and has not yet reached a steady state.  We have plotted the time evolution of kinetic energy, charge drops and outflux in Figure $\ref{fig:QuantitiesBifurcated}$ for $Re = 600.$

\begin{figure}[ht!]
    \centering
    \includegraphics[width=1\linewidth]{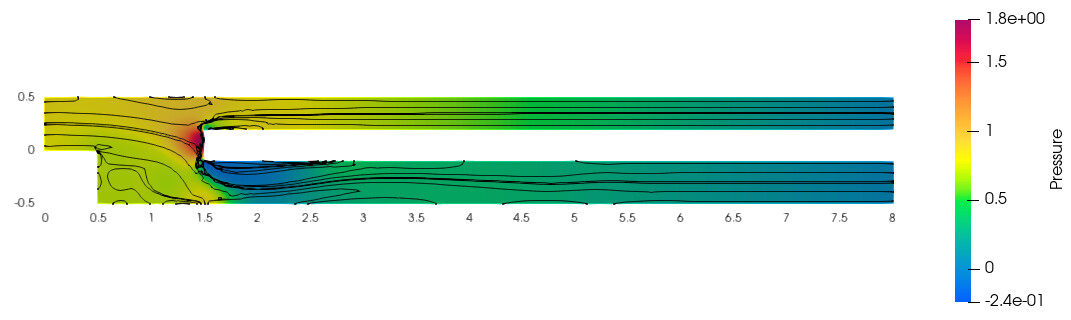}
    \caption{Example \ref{Bifur}: Pressure and streamlines of vorticity for $Re = 600$ at final time $T = 7s.$}
    \label{fig:PressureIsoBifurcated}
\end{figure}

Moreover, we show in Figure \ref{fig:PressureIsoBifurcated}  the pressure and streamlines of the vorticity for $Re = 600$ at final time $T = 7s.$

\begin{figure}[ht!]
\centering
\begin{subfigure}[b]{0.45\textwidth}
\includegraphics[width=\textwidth]{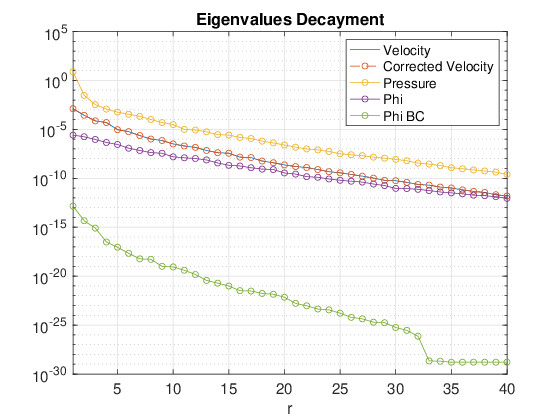}
\end{subfigure}
\begin{subfigure}[b]{0.45\textwidth}
\includegraphics[width=\textwidth]{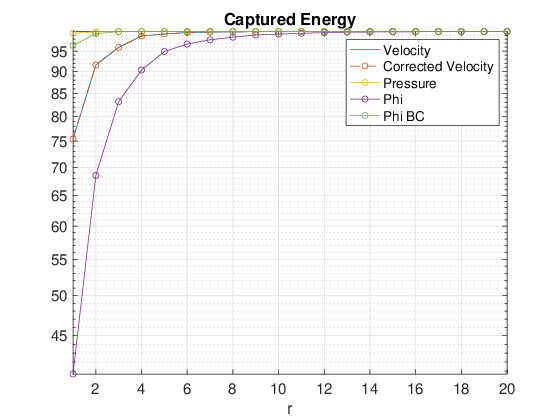}
\end{subfigure}
\caption{Example \ref{Bifur}: Decay of the normalized POD eigenvalues (left) and percent ratio of captured energy (right).}
\label{fig:eigenvaluesBifurcated}
\end{figure}

We consider, as mentioned before, the time and the Reynolds number as parameters to build the ROM. We take $5$ samples of a uniform partition of Reynolds number range $\mathcal{D} = [400,600]$. Corcerning the time parameter, we collect $201$ snapshots in the time interval $[5,7]s,$ in which the flow has not yet reached the steady state. Finally, we obtain $M = 5 \times 201 = 1005$ snapshots for each unknown field. The POD modes are generated in  $L^2-$norm for velocity field and $H^1-$norm for pressure field. Figure \ref{fig:eigenvaluesBifurcated} shows the eigenvalues decay (left) associated to the correlation matrix of the velocity, predicted velocity, pressure, pressure correction and pressure boundary condition. Moreover, we can also see in Figure \ref{fig:eigenvaluesBifurcated} the percentaje of captured energy (right) computed as $100 (\sum_{k=1}^r \lambda_k )/ (\sum_{k=1}^E \lambda_k),$ where $\lambda_k$ are the corresponding eigenvalues and $E$ the rank of the corresponding data set. We can see in Figure \ref{fig:eigenvaluesBifurcated} that with just $10$ modes for each field we obtain more than $95\%$ of the captured energy.

\subsubsection{Numerical results}
Once the POD modes are generated, we study the performances of the POD-ROM and the hybrid POD-ROM, using the standard Radial Basis Method with time and Reynolds number as parameter explained in \ref{RBFsection}, in terms of the accuracy and efficiency.

\begin{figure}[ht!]
\centering
\begin{subfigure}[b]{0.45\textwidth}
\includegraphics[width=\textwidth]{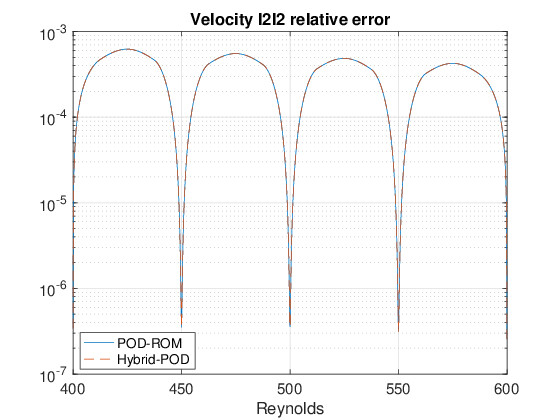}
\end{subfigure}
\begin{subfigure}[b]{0.45\textwidth}
\includegraphics[width=\textwidth]{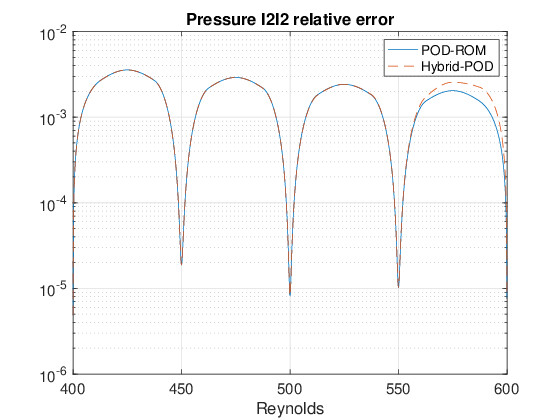}
\end{subfigure}
\caption{Example \ref{Bifur}: $l^2(L^2)$ relative error for velocity (left) and pressure (right) in the time interval $[5,7]s$ for different values of Reynolds number inside $\mathcal{D}.$}
\label{fig:BifurcatedReynoldsROM}
\end{figure}

In particular, we show in Figure \ref{fig:BifurcatedReynoldsROM} the $l^2(L^2)$ relative error for velocity and pressure in the time interval $[5,7]s$ for different values of Reynolds numbers inside the range $\mathcal{D} = [400,600],$ but different from the chosen samples. We can note from Figure \ref{fig:BifurcatedReynoldsROM} that we almost obtain the same error level in velocity and pressure for both ROMs.

\begin{figure}[ht!]
\centering
\begin{subfigure}[b]{0.45\textwidth}
\includegraphics[width=\textwidth]{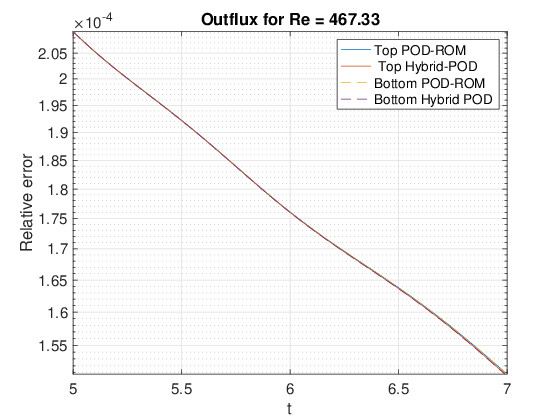}
\end{subfigure}
\begin{subfigure}[b]{0.45\textwidth}
\includegraphics[width=\textwidth]{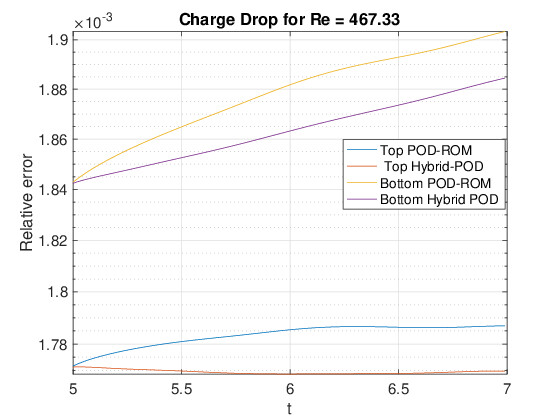}
\end{subfigure}
\begin{subfigure}[b]{0.45\textwidth}
\includegraphics[width=\textwidth]{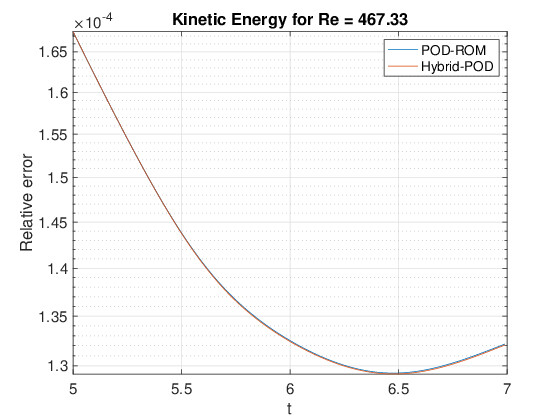}
\end{subfigure}
\caption{Example \ref{Bifur}: Temporal evolution of relative error between FOM and ROM of the outflux (top left figure), charge drops (top right figure) and kinetic energy (bottom figure) for $Re = 467.33$}
\label{fig:Re467}
\end{figure}

\begin{figure}[ht!]
\centering
\begin{subfigure}[b]{0.45\textwidth}
\includegraphics[width=\textwidth]{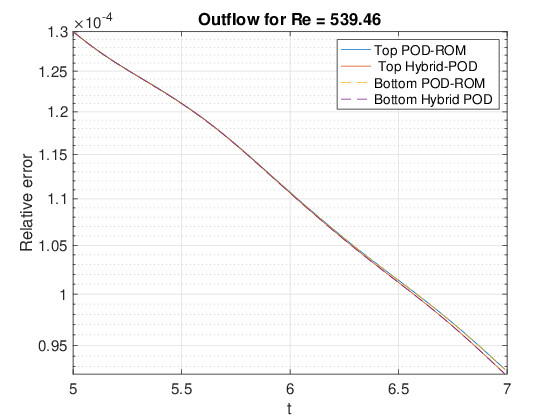}
\end{subfigure}
\begin{subfigure}[b]{0.45\textwidth}
\includegraphics[width=\textwidth]{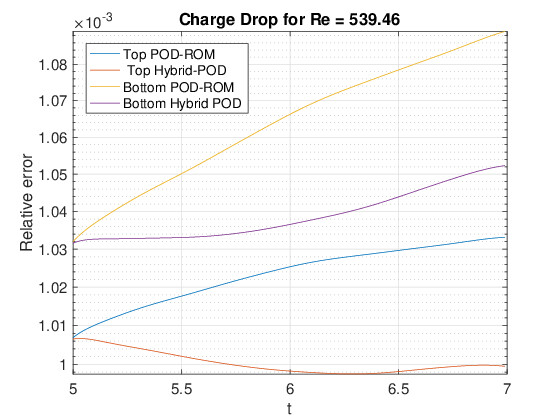}
\end{subfigure}
\begin{subfigure}[b]{0.45\textwidth}
\includegraphics[width=\textwidth]{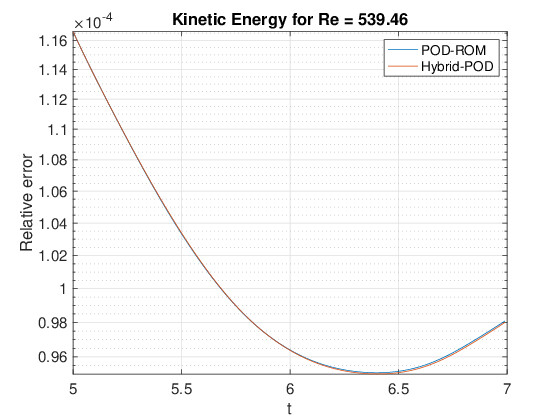}
\end{subfigure}
\caption{Example \ref{Bifur}: Temporal evolution of relative error between FOM and ROM of the outflux (upper left figure), charge drops (upper right figure) and kinetic energy (down figure) for $Re = 539.46$}
\label{fig:Re539}
\end{figure}

However, for a more comprehensive evaluation of both ROMs, Figures \ref{fig:Re467} and \ref{fig:Re539}  present the temporal evolution of the relative error between FOM and ROM for quantities of interest such as kinetic energy, charge drops, and outflux, considering two Reynolds number values, \( Re = 466 \) and \( Re = 539 \), which lie within the trial interval, but differ from the selected samples.

We can see in Figures \ref{fig:Re467} and \ref{fig:Re539} that both ROMs almost get the same levels of accuracy in terms of relative error, obtaining in both cases errors below $2 \cdot 10^{-4}$ for outflux and kinetic energy and $2 \cdot 10^{-3}$ for charge drops.

Finally, in Table \ref{tab:CPUtimeBifurcatedROM}  we show the CPU time of FOM and ROMs in the time interval $[5,7]s$ and its corresponding speed up. Note that the hybrid model is slightly faster than the fully intrusive one while almost achieving the same level of accuracy.

\begin{table}[ht!]
\centering
\begin{tabular}{|l|l|l|l|l|}
\hline
CPU time (s)      & $Re = 467.33$ & Speed Up & $Re = 539.46$ & Speed Up \\ \hline
Full Order Model  & $346.96s$     & $1$      & $344.38$      & 1        \\ \hline
Intrusive POD-ROM & $0.214s$      & $1621$   & $0.227s$      & $1517$   \\ \hline
Hybrid POD-ROM    & $0.201s$      & $1726$   & $0.198s$      & $1739$   \\ \hline
\end{tabular}
\caption{CPU time of FOM and ROMs in the time interval $[5,7]s$ for $Re = 467.33$ and $539.46$, and its corresponding speed up.}
\label{tab:CPUtimeBifurcatedROM}
\end{table}

\subsection{Flow past a cylinder }\label{FlowCylinder}
In the next example, we consider the benchmark problem of 2D unsteady flow around a cylinder with circular cross-section \cite{SchaferTurek96}, and we take the same setup for numerical simulations as in \cite{samujulia}. 

\subsubsection{Checking of outflow boundary position and POD modes}
Our aim is to compare the performance of the fully-discrete time-splitting FOM described in Section \ref{FETimeSplit} with both ROM counterparts described in Section \ref{ROMtimesplit} and \ref{HybridROM}. To analyze the effect of the position of the outflow boundary, we consider three computational domains where the distance between the inlet and the outlet is respectively $2.2$ (full domain), $1.1$ (half domain), and $0.6$ (small domain), as we can see in Figure \ref{fig:Meshs}:

\begin{figure}[ht!]
\centering
\begin{subfigure}[b]{0.7\textwidth}
\includegraphics[width=\textwidth]{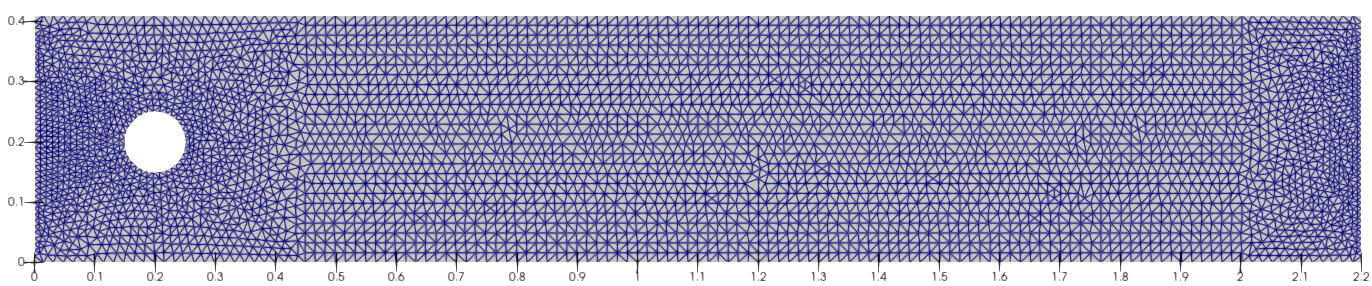}
\end{subfigure}
\begin{subfigure}[b]{0.7\textwidth}
\includegraphics[width=\textwidth]{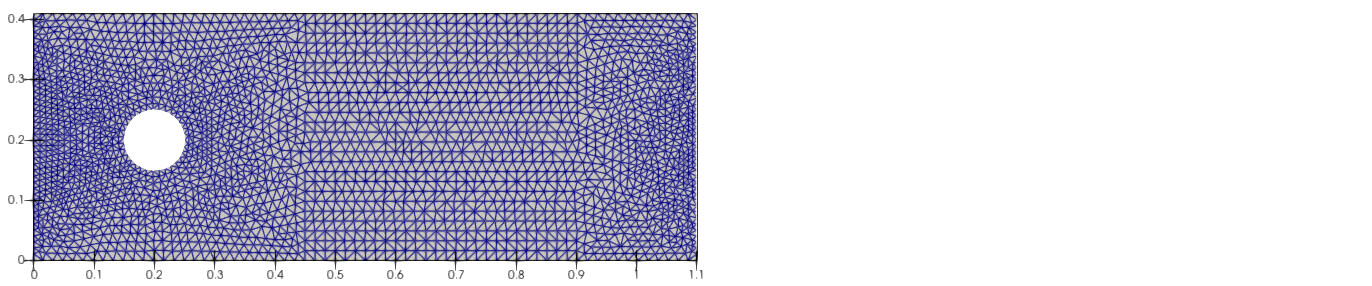}
\end{subfigure}
\begin{subfigure}[b]{0.7\textwidth}
\includegraphics[width=\textwidth]{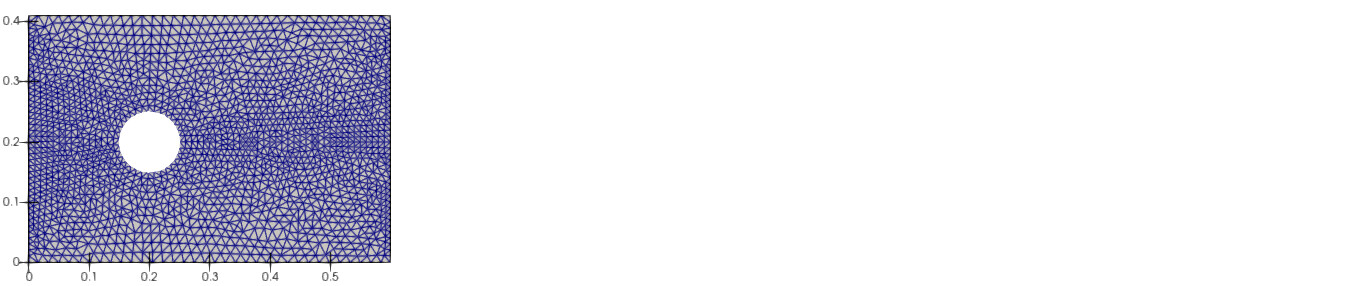}
\end{subfigure}
\caption{Example \ref{FlowCylinder}: Three different computational grids depending on the distance between inlet and outlet.}
\label{fig:Meshs}
\end{figure}

The FOM is the one described in Section \ref{FETimeSplit}, with a spatial discretization using  $\mathbb{P}^2$ FE for the velocity and predicted velocity and $\mathbb{P}^1$ FE for the pressure, pressure-continuity correction and pressure-continuity boundary correction on a relatively coarse computational grid (see Figure \ref{fig:Meshs}), for which  the maximum grid sizes are $h_1 = 2.6 \cdot 10^{-2}m, h_2 = 2.7 \cdot 10^{-2}m,$ and $ h_3 = 2.4 \cdot 10^{-2}m$, respectively for the full, half and small domain (from top to bottom in Figure \eqref{fig:Meshs}), resulting in $33.518, 19.950$ and $14.626$  degrees of freedom for velocities and $4.284, 2.555$ and $1.874$ for pressures.

\begin{figure}[ht!]
\centering
\begin{subfigure}[b]{0.45\textwidth}
\includegraphics[width=\textwidth]{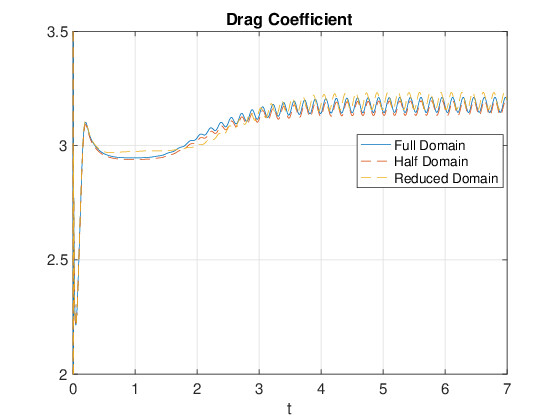}
\end{subfigure}
\begin{subfigure}[b]{0.45\textwidth}
\includegraphics[width=\textwidth]{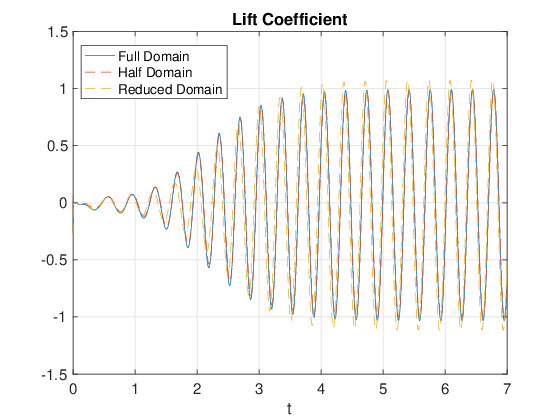}
\end{subfigure}
\begin{subfigure}[b]{0.45\textwidth}
\includegraphics[width=\textwidth]{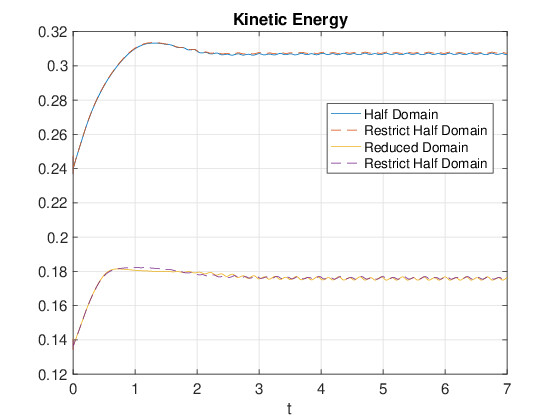}
\end{subfigure}
\caption{ Example \ref{FlowCylinder}: Temporal evolution of drag (top left figure) and lift (top right figure) coefficients for the FOM solution computed on the three computational domains and temporal evolution of kinetic energy (bottom figure) of the full domain restricted into the other computational domains.}
\label{fig:QuantitiesFOM}
\end{figure}

In the FOM simulations, an impulsive start is performed, i.e the initial conditions are a zero velocity and zero pressure fields, and the time step considered is $\Delta t = 2 \cdot 10^{-3}s.$ Time integration is performed till a final time $T = 7s.$ We present in Figure \ref{fig:QuantitiesFOM} the time evolution of the drag and lift coefficients obtained from the full-order model (FOM) simulations across the three computational domains. In addition, the figure also displays the evolution of kinetic energy, as well as the lift and drag coefficients restricted to the respective subdomains, over the time interval $[5,7]s$. This allows for a detailed comparison of the flow behavior and physical quantities across different domain configurations depending on the distance between inlet and the outlet.

\begin{figure}[ht!]
\centering
\begin{subfigure}[b]{0.7\textwidth}
\includegraphics[width=\textwidth]{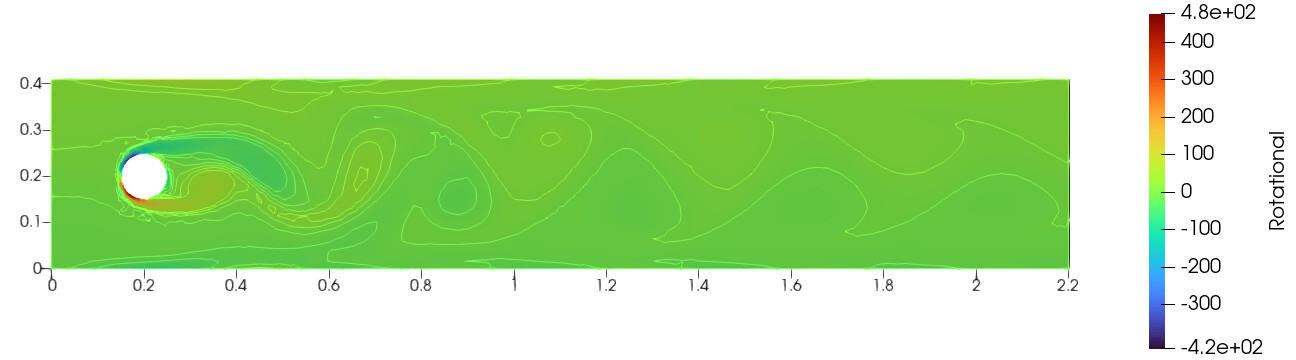}
\end{subfigure}
\begin{subfigure}[b]{0.7\textwidth}
\includegraphics[width=\textwidth]{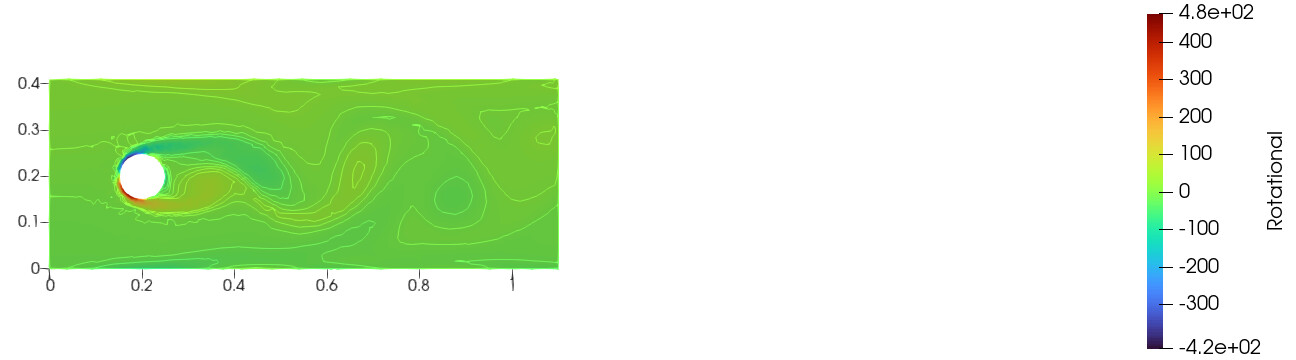}
\end{subfigure}
\begin{subfigure}[b]{0.7\textwidth}
\includegraphics[width=\textwidth]{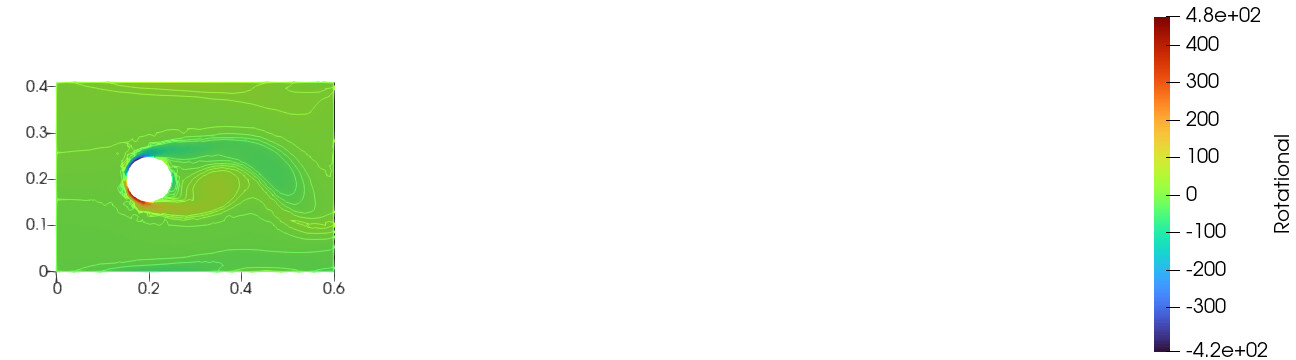}
\end{subfigure}
\caption{Example \ref{FlowCylinder}: Vorticity magnitude and contours for the FOM solution computed on the three computational domains at final simulation time $T=7s.$}
\label{fig:RotationalFOM}
\end{figure}

In Figure \ref{fig:RotationalFOM}, we show the vorticity magnitude and contours for the FOM solution computed on the three computational domains at final simulation time $T = 7s.$ We note that the position of the outflow boundary condition does not induce distortion of the flow. Hence, thanks to the non-standard treatment of the open boundary condition within the fully discrete time-splitting FOM, it would be sufficient just to compute the solution in the reduced domain in order to obtain quantities of interest such as drag and lift coefficient, with considerable saving in CPU time. Indeed, we show in Table \ref{tab:CPUtimeFOM} the CPU time of the FOM time-splitting method on the three computational domains in the time interval $[0,7]s.$ We observe a gain of over $38\%$ CPU time using the half domain and of nearly $54\%$ using the small domain. 

\begin{table}[ht!]
\centering
\begin{tabular}{|l|lll|lll|}
\hline
CPU time (s)   & \multicolumn{3}{l|}{}           & \multicolumn{3}{l|}{Speed Up} \\ \hline
Full Domain    & \multicolumn{3}{l|}{$5027.21s$} & \multicolumn{3}{l|}{$1$}      \\ \hline
Half Domain    & \multicolumn{3}{l|}{$3092.97s$} & \multicolumn{3}{l|}{$1.625$}  \\ \hline
Reduced Domain & \multicolumn{3}{l|}{$2313.57s$} & \multicolumn{3}{l|}{$2.172$}  \\ \hline
\end{tabular}
\caption{Example \ref{FlowCylinder}: CPU time of fully discrete time-splitting method on the three computational domains in the time interval $[0,7]s$ and its corresponding speed up.}
\label{tab:CPUtimeFOM}
\end{table}

\begin{figure}[ht!]
\centering
\begin{subfigure}[b]{0.45\textwidth}
\includegraphics[width=\textwidth]{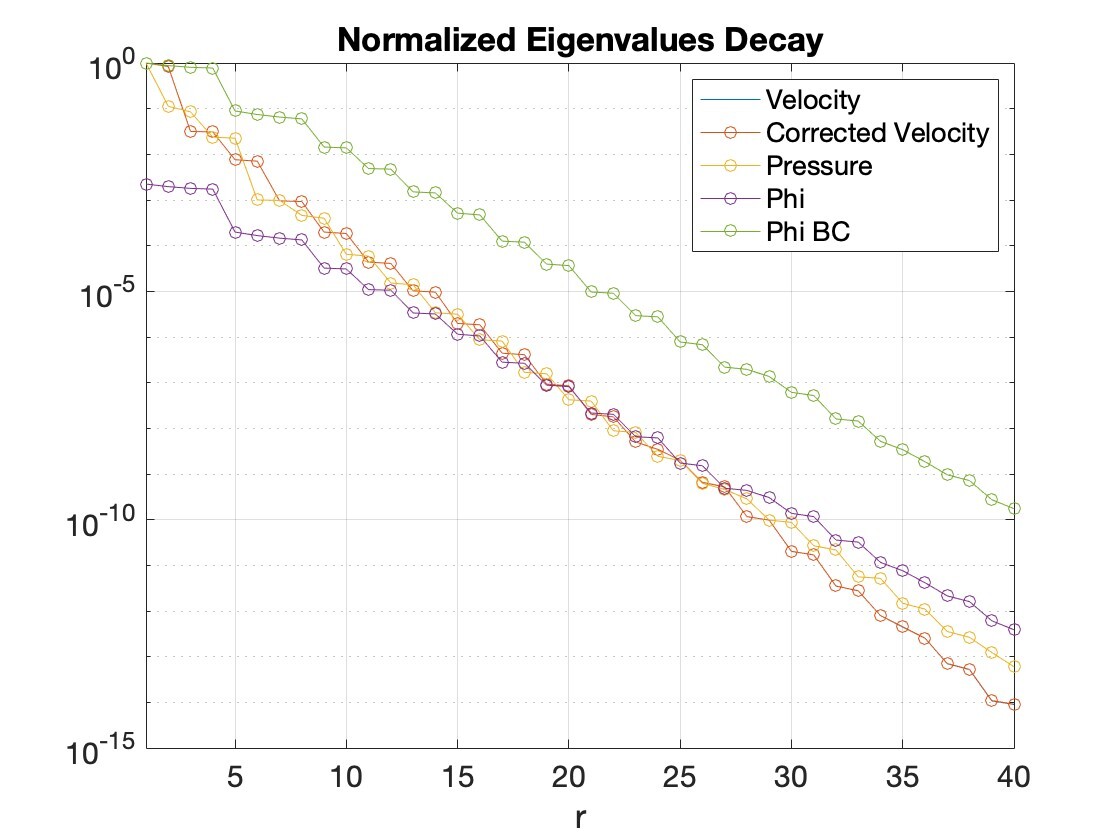}
\end{subfigure}
\begin{subfigure}[b]{0.45\textwidth}
\includegraphics[width=\textwidth]{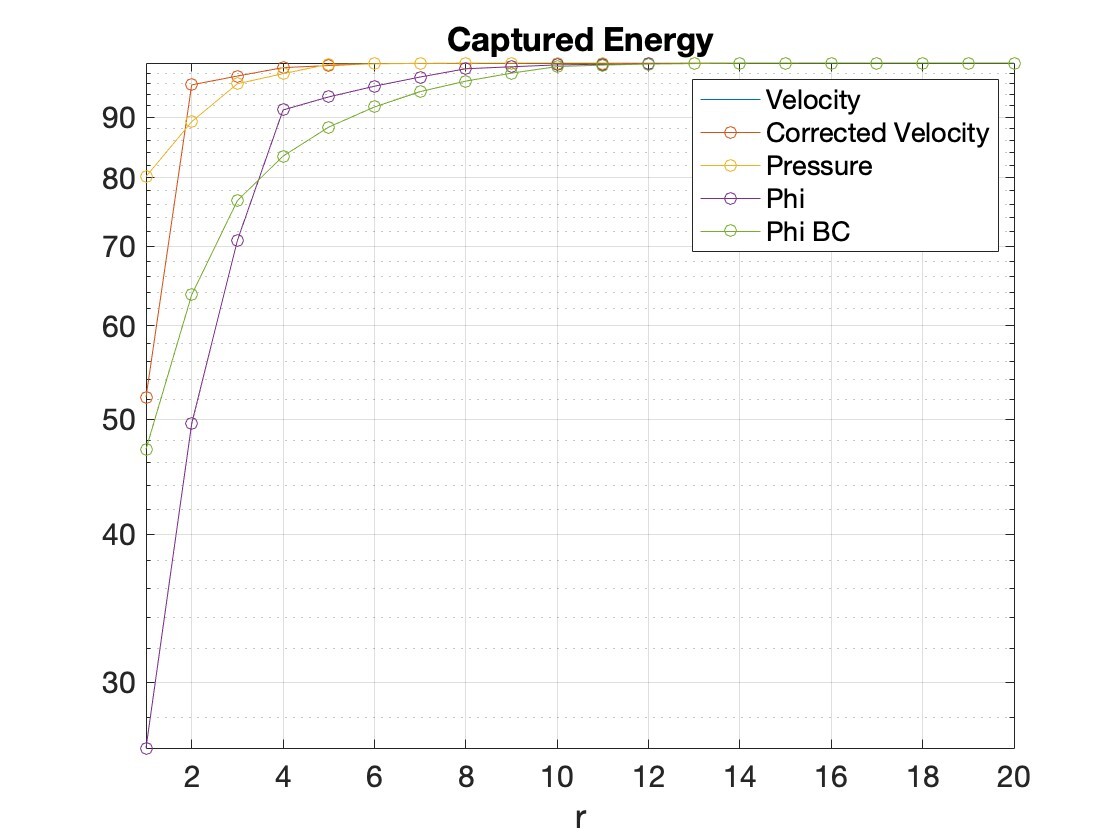}
\end{subfigure}
\caption{Example \ref{FlowCylinder}: Decay of the normalized POD eigenvalues (left) and percent ratio of captured energy (right).}
\label{fig:eigenvalues}
\end{figure}

The ROM is applied just to the reduced domain. The POD modes are generated in $L^2-$norm for the velocity field and $H^1-$norm for the pressure field by the method detailed in Section \ref{OfflinePhase}, by storing every FOM solution from $t = 5s$ up to complete one period. The full period length of the statistically steady state is $0.332s$ for $Re = 100,$ so we collect $167$ snapshots for each unknown. In Figure \ref{fig:eigenvalues}, we can see the decay of the normalized POD eigenvalues (left) associated to the correlation matrix of the velocity, predicted velocity, pressure, pressure correction and pressure boundary condition. Figure \ref{fig:eigenvalues} also displays the corresponding percentaje of captured energy (right) computed by $100 (\sum_{k=1}^r \lambda_k) / (\sum_{k=1}^E \lambda_k),$ where $\lambda_k$ are the corresponding eigenvalues and $E$ the rank of the corresponding data set. Note that with 8 POD modes we already capture more than $95\%$ of the energy for each field. 

Once the POD modes are generated, the POD-ROM and Hybrid-ROM, using RBF extrapolation explained in \ref{RBFextrapolation} (in which the parameters are the POD coefficients for the velocity $\bm{\tilde{u}}_r$), for the fully discrete time-splitting FOM are constructed. Both ROMs are run from $t=5s$ up to 6 periods, i.e. in the time interval $[5,6.992]s$, using 20 POD modes for each field.

\subsubsection{Numerical results on the reduced domain}
As we mentioned before, the POD-ROM and Hybrid-ROM, using RBF extrapolation, for the fully discrete time-splitting FOM are applied just to the reduced domain and tested in the time interval $[5,6.992]s$. 

\begin{figure}[ht!]
\centering
\begin{subfigure}[b]{0.45\textwidth}
\includegraphics[width=\textwidth]{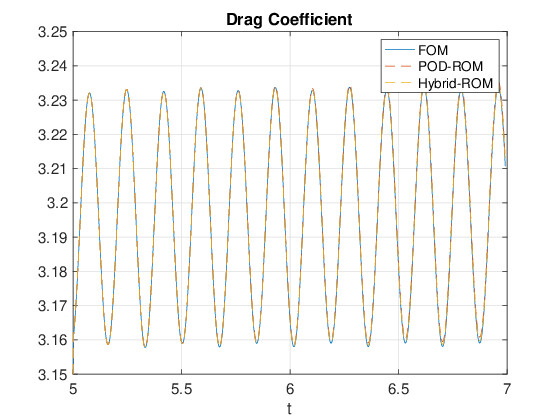}
\end{subfigure}
\begin{subfigure}[b]{0.45\textwidth}
\includegraphics[width=\textwidth]{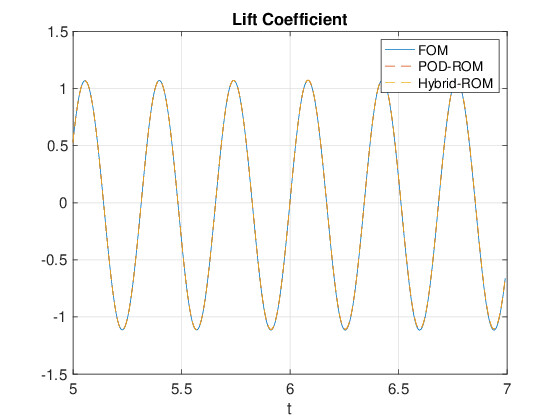}
\end{subfigure}
\begin{subfigure}[b]{0.45\textwidth}
\includegraphics[width=\textwidth]{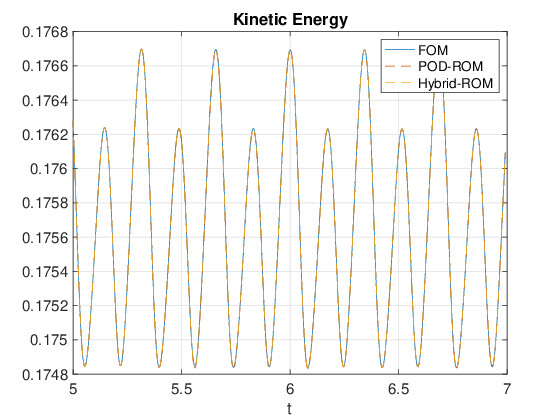}
\end{subfigure}
\caption{Example \ref{FlowCylinder}: Time evolution of the drag coefficient (top left figure), lift coefficient (top right figure) and kinetic energy (bottom figure) for FOM and POD-ROM.}
\label{fig:QuantitiesROM}
\end{figure}

We show the numerical results for kinetic energy, drag and lift coefficient in Figure \ref{fig:QuantitiesROM}. We can observe that both ROMs accurately capture the three quantities of interest.   

\begin{figure}[ht!]
\centering
\begin{subfigure}[b]{0.45\textwidth}
\includegraphics[width=\textwidth]{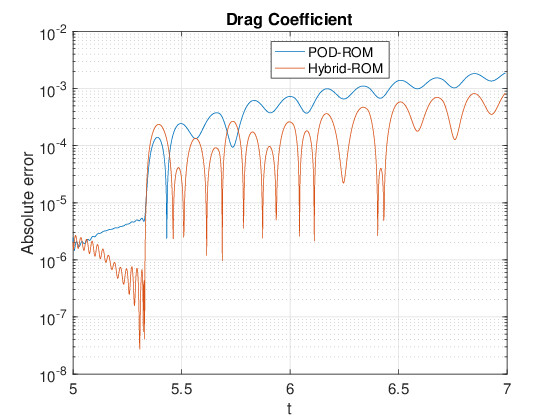}
\end{subfigure}
\begin{subfigure}[b]{0.45\textwidth}
\includegraphics[width=\textwidth]{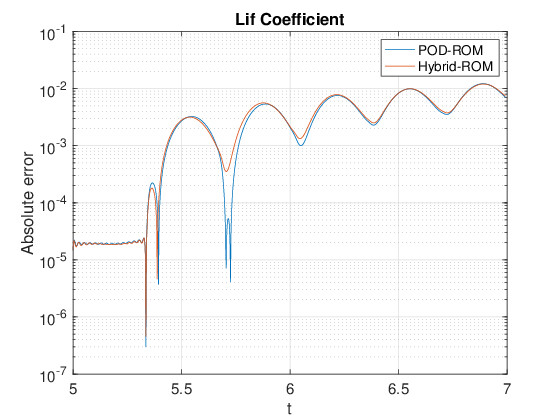}
\end{subfigure}
\begin{subfigure}[b]{0.45\textwidth}
\includegraphics[width=\textwidth]{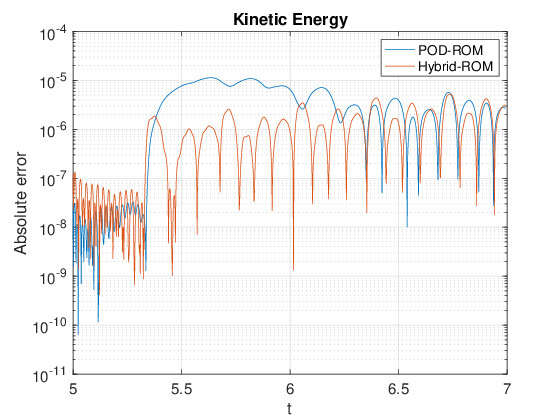}
\end{subfigure}
\caption{Example \ref{FlowCylinder}: Time evolution of the absolute error of drag coefficient (top left figure), lift coefficient (top right figure) and kinetic energy (bottom figure) for FOM and POD-ROM.}
\label{fig:ErrorQuantitiesROM}
\end{figure}

To better quantify these results, Figure \eqref{fig:ErrorQuantitiesROM} shows the absolute error of these quantities of interest by comparing the ROM with the FOM. Both ROMs achieve quite close error levels for drag and lift coefficients. However, for kinetic energy, we can observe that the fully intrusive approach guarantee better error levels. Moreover, there is an increase in the error once the temporal extrapolation begins. However, the error remains small, on the order of $10^{-3}.$

\begin{figure}[ht!]
\centering
\begin{subfigure}[b]{0.45\textwidth}
\includegraphics[width=\textwidth]{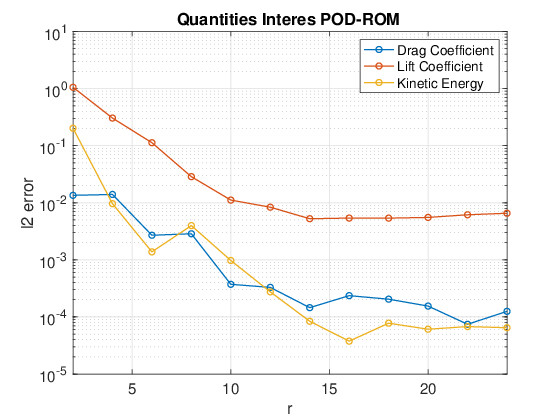}
\end{subfigure}
\begin{subfigure}[b]{0.45\textwidth}
\includegraphics[width=\textwidth]{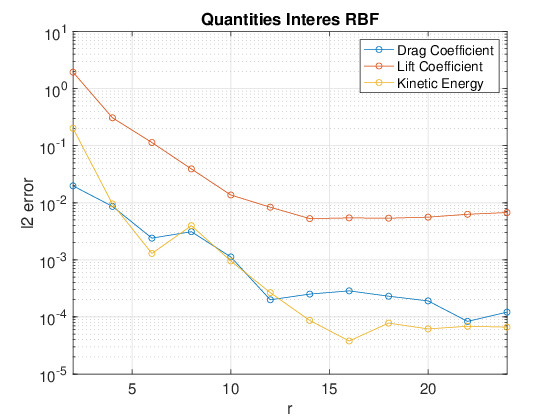}
\end{subfigure}
\caption{Example \ref{FlowCylinder}: The $L^2$ relative errors for the kinetic energy, drag and lift coefficients versus number of modes $r$. The error is computed between the FOM quantities and ROM quantities for the time interval $[5,6.992]s.$}
\label{fig:l2Quantities}
\end{figure}

Furthermore, we show in Figure \ref{fig:l2Quantities} the $L^2$ in time relative errors in the time interval $[5,6.992]s$ depending on the number of POD modes $r.$ The $L^2$ relative error is defined as follows: 
\begin{equation}
    \epsilon_{L^2} = \dfrac{\|\mu(t) - \mu^*(t) \|_{L^2(T_1,T_2)}}{\|\mu(t)\|_{L^2(T_1,T_2)}},
\end{equation}
where $\mu(t)$ is a quantity of interest computed with the FOM and $\mu^*$ is the corresponding quantity computed with the ROM. In particular, looking at Figure \ref{fig:l2Quantities}, for both ROMs we can observe that drag coefficient and kinetic energy errors are around $0.01\%$. However, the discrepancy in the lift coefficient remains below $1\%$, which is likely due to a slight phase error between the FOM and the POD-ROM during the extrapolation interval.

\begin{figure}[ht!]
\centering
\begin{subfigure}[b]{0.45\textwidth}
\includegraphics[width=\textwidth]{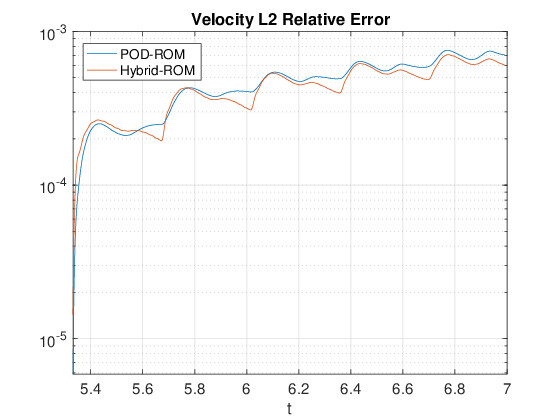}
\end{subfigure}
\begin{subfigure}[b]{0.45\textwidth}
\includegraphics[width=\textwidth]{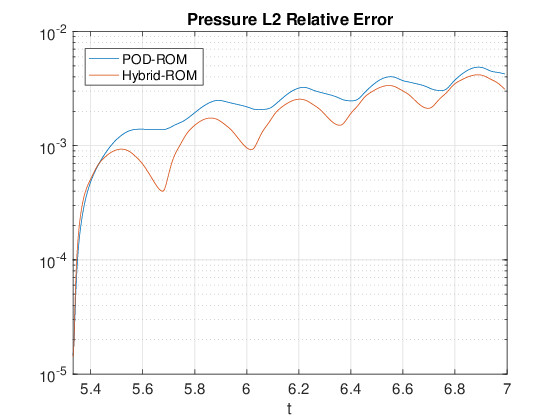}
\end{subfigure}
\begin{subfigure}[b]{0.45\textwidth}
\includegraphics[width=\textwidth]{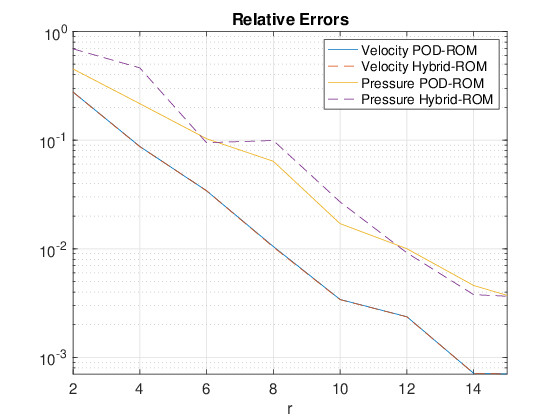}
\end{subfigure}
\caption{Example \ref{FlowCylinder}: Temporal evolution of relative errors (in $L^2-$norm) for velocity and pressure (top figures) and $l^2(L^2)$ relative errors of velocity and pressure depending on the number of POD basis $r$ (bottom figure).}
\label{fig:REerror}
\end{figure}

To better assess the numerical accuracy of the proposed ROMs, we show in Figure \ref{fig:REerror} the relative error between FOM and ROM for velocity and pressure in the extrapolation time interval (above) and
$l^2(L^2)$ relative errors depending on the number of POD modes $r.$ 
As we can see in Figure \ref{fig:REerror} (left), the relative errors for velocity and pressure take values close to $10^{-3}$ and $10^{-2}$, respectively.  On the other hand, we observe that the error decay stagnates when the number of modes exceeds 14. Let us recall that this is a common phenomenon in POD-Galerkin ROMs.

We show in Table \ref{tab:CPUtimeROM} the CPU time of the fully discrete time-splitting method and the reduced order model on the reduced domain in the time interval $[5,6.992]s.$ 

\begin{table}[ht!]
\centering
\begin{tabular}{|l|lll|lll|}
\hline
CPU time (s)        & \multicolumn{3}{l|}{}           & \multicolumn{3}{l|}{Speed Up} \\ \hline
Full Order Model    & \multicolumn{3}{l|}{$661.02s$} & \multicolumn{3}{l|}{$1$}      \\ \hline
Intrusive POD-ROM & \multicolumn{3}{l|}{$0.427s$}   & \multicolumn{3}{l|}{$1547.80$} \\ \hline
Hybrid POD-ROM & \multicolumn{3}{l|}{$0.93218s$}   & \multicolumn{3}{l|}{$709.11$} \\ \hline
\end{tabular}
\caption{Example \ref{FlowCylinder}: CPU time of FOM and ROMs on the reduced computational domain in the time interval $[5,6.992]s$ and its corresponding speed up.}
\label{tab:CPUtimeROM}
\end{table}

\begin{figure}[ht!]
\begin{subfigure}[t]{0.35\textwidth}
    \includegraphics[width=\linewidth]{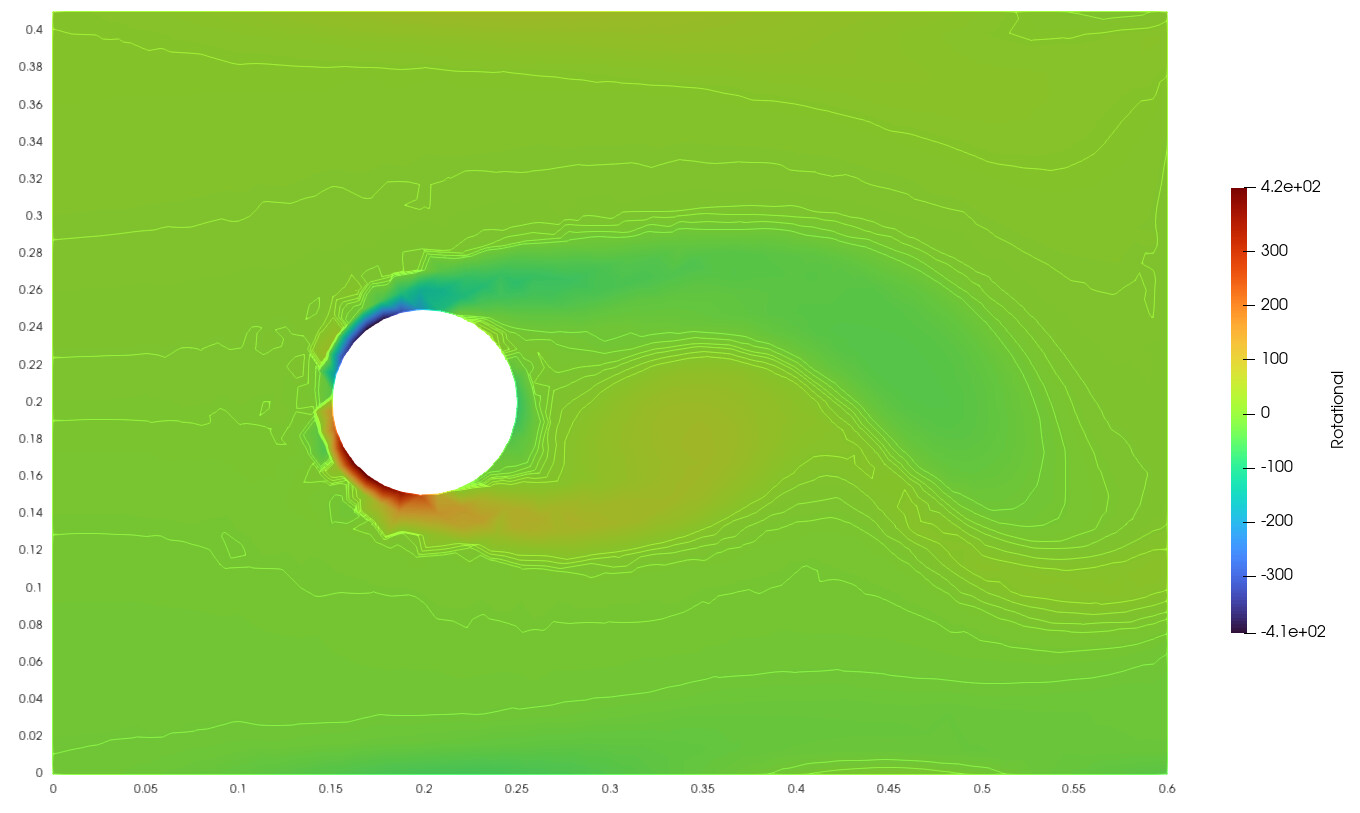}
\end{subfigure}
\begin{subfigure}[t]{0.35\textwidth}
    \includegraphics[width=\linewidth]{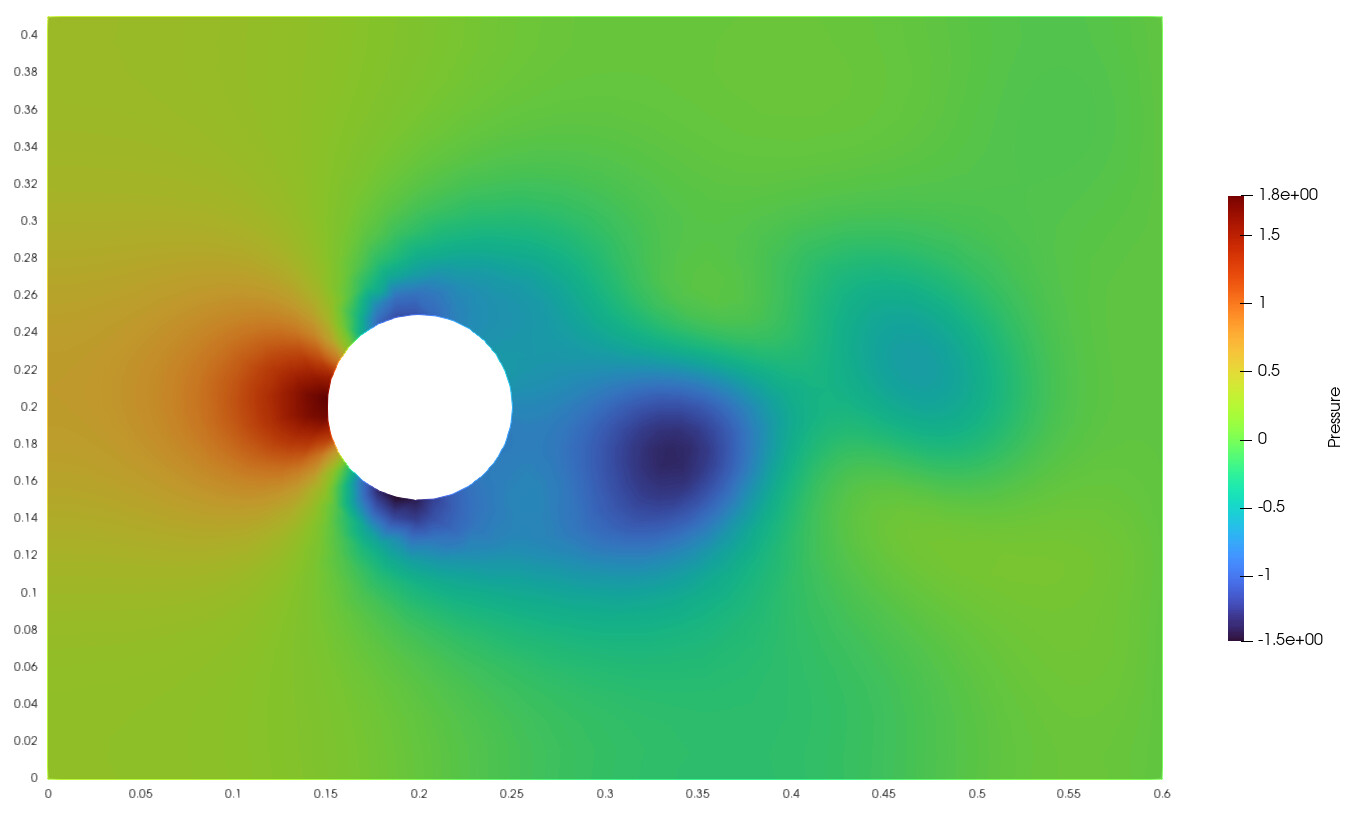}
\end{subfigure}
\begin{subfigure}[t]{0.35\textwidth}
    \includegraphics[width=\linewidth]{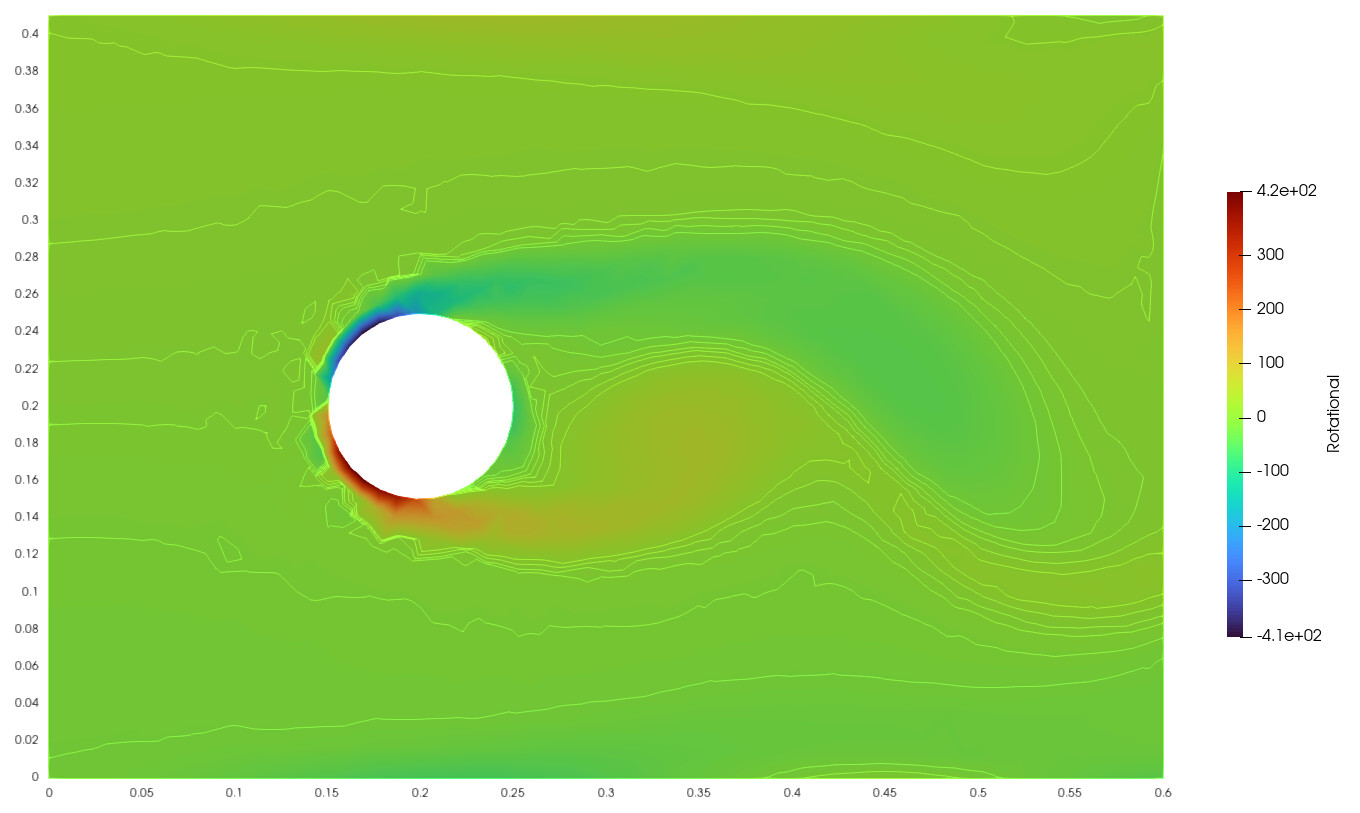}
\end{subfigure}\hfill
\begin{subfigure}[t]{0.35\textwidth}
    \includegraphics[width=\linewidth]{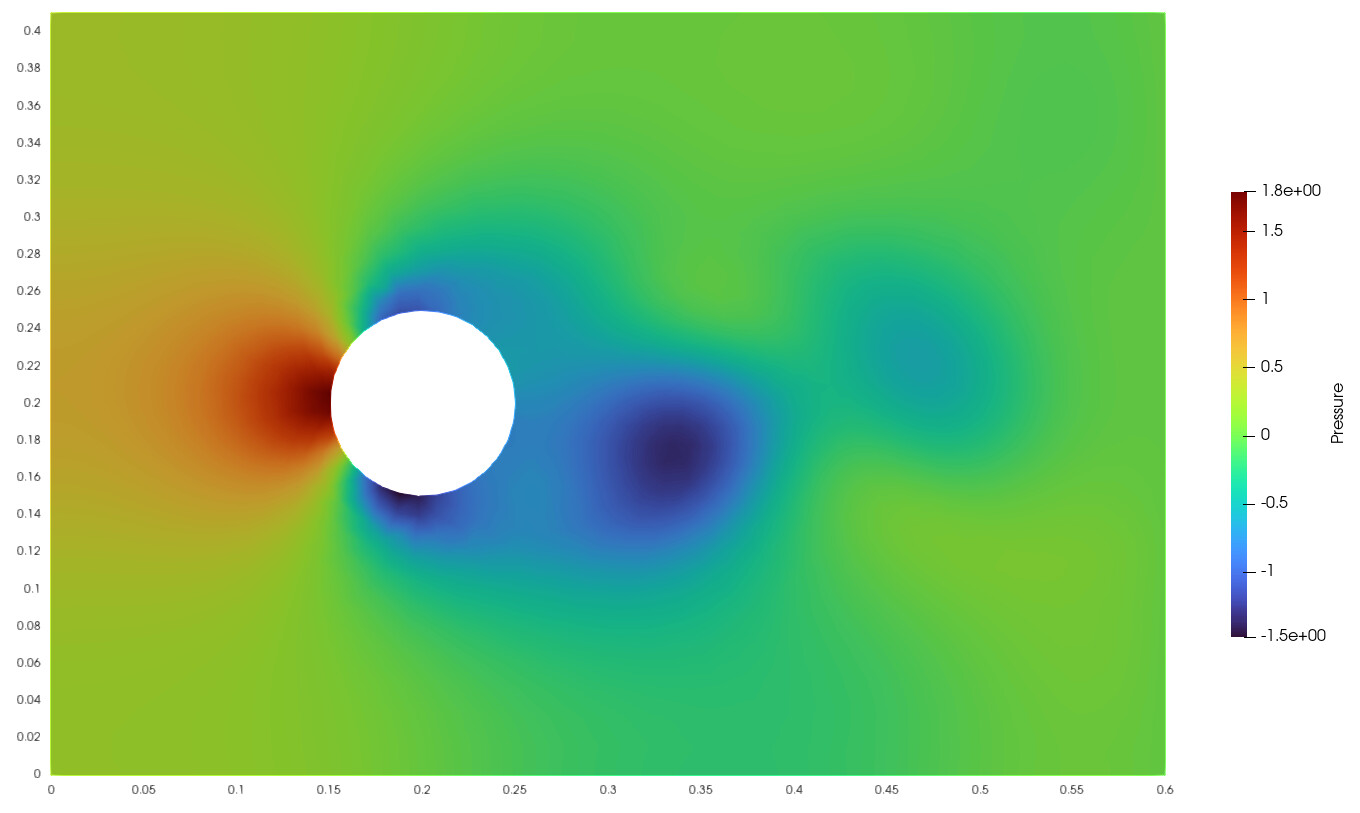}
\end{subfigure}
\begin{subfigure}[t]{0.35\textwidth}
    \includegraphics[width=\linewidth]{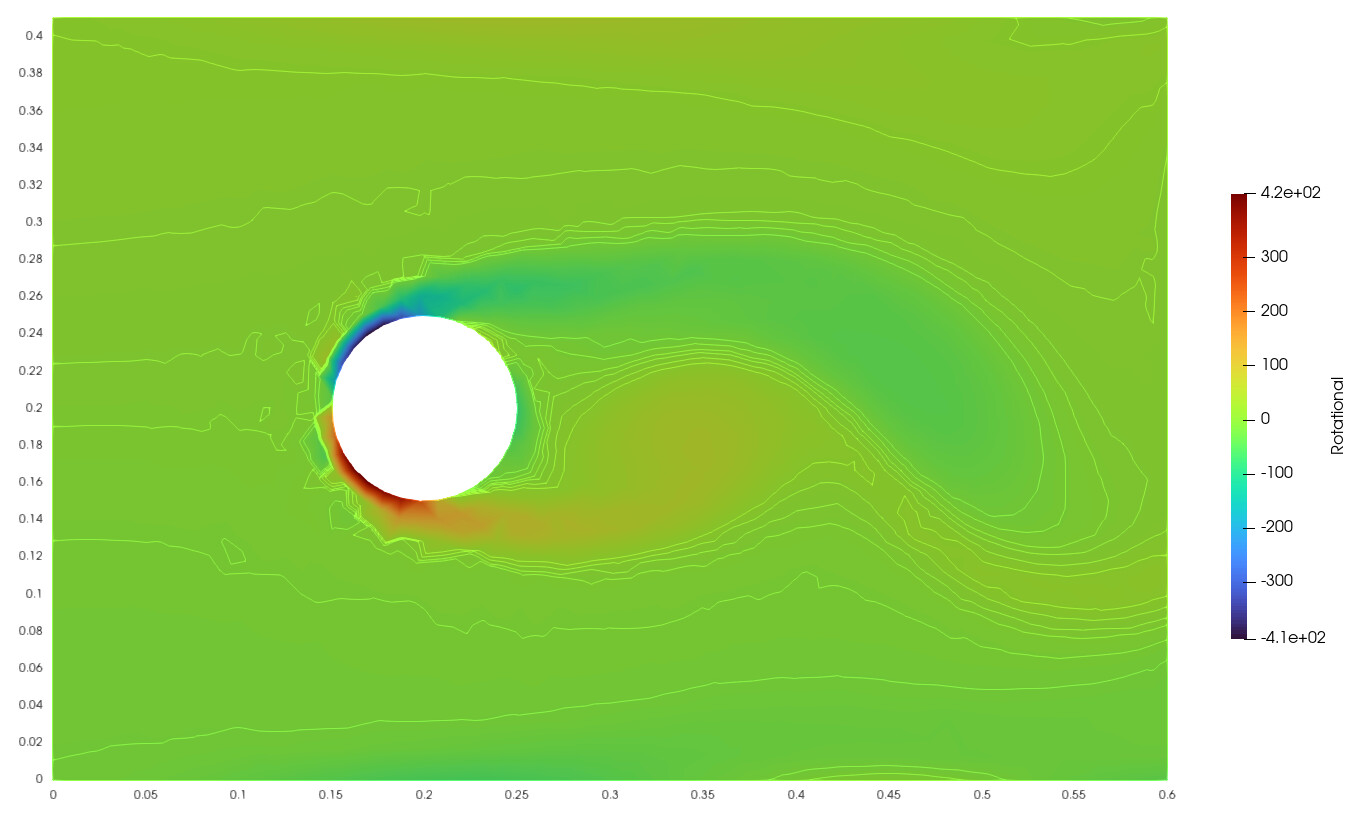}
\end{subfigure}\hfill
\begin{subfigure}[t]{0.35\textwidth}
    \includegraphics[width=\linewidth]{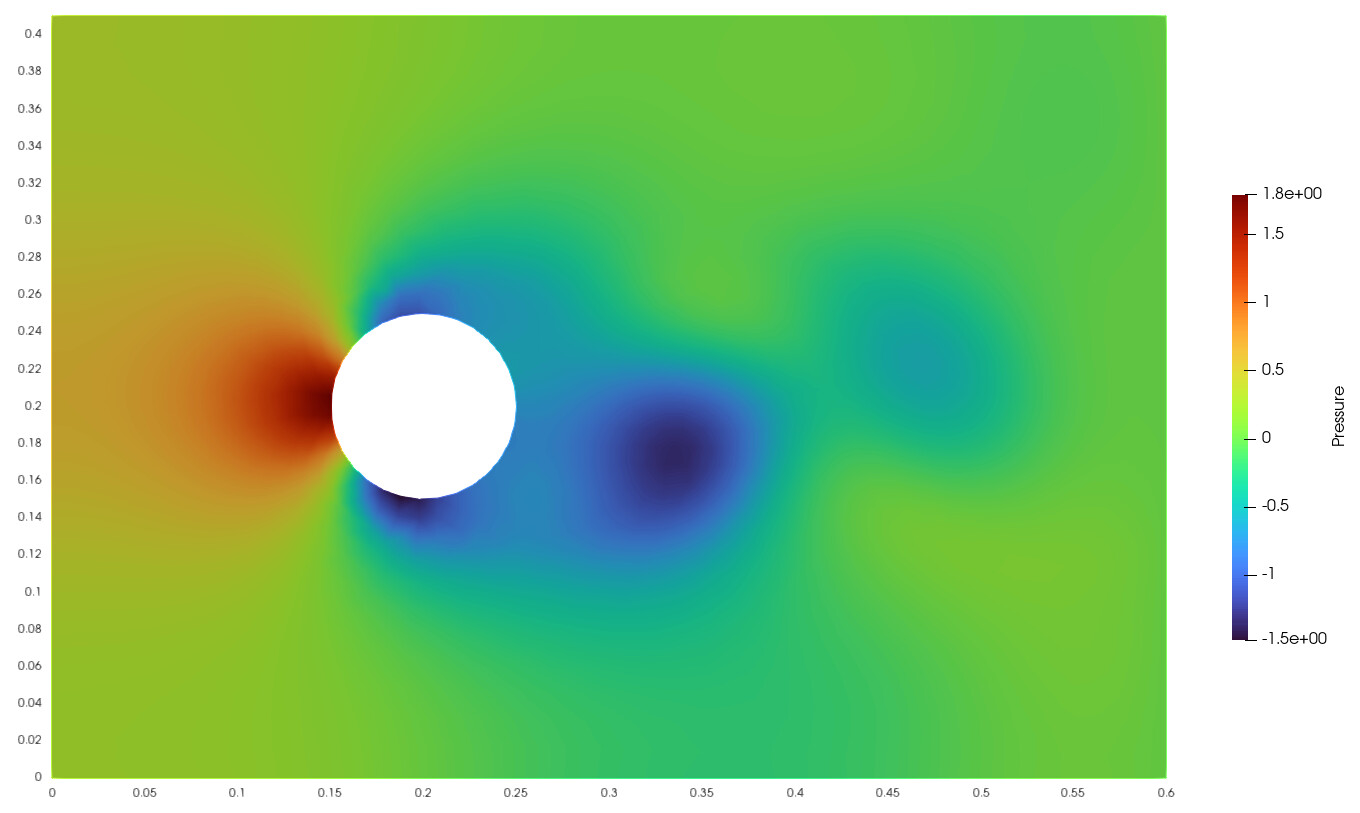}
\end{subfigure}
\caption{Example \ref{FlowCylinder}: Vorticity magnitude and contours (left) and the pressure field (right) of FOM (top), POD-ROM (medium) and Hybrid POD-ROM (bottom) after the extrapolation of six periods over the reduced computational domain.}
\label{fig:RotationalROM}
\end{figure}

Finally, we show in Figure \ref{fig:RotationalROM} the vorticity magnitude and contours (left) and the pressure field (right) of FOM (top), POD-ROM (medium) and Hybrid POD-ROM (bottom) after the extrapolation of six periods over the reduced computational domain. The results for both ROMs (almost identical) show that using a reduced computational domain does not induce distortion of the vortices nor disturbs the flow around the cylinder. 

\section{Conclusions}\label{Conclusion}

In this work, we present a POD-ROM for the time-splitting method applied to the Navier-Stokes equations with open boundary conditions. In the proposed time-splitting scheme, the pressure-correction method is suitable for handling open boundary conditions without compromising the accuracy the solution, when a reduced computational domain is considered. The combination of this scheme with a POD-ROM enhances to exploit efficiency of the simulations maintaing accuracy.

Additionally, we propose two reduced-order modeling approaches: the standard POD-ROM projection-based method (fully intrusive) and the hybrid POD-ROM, which is partially projection-based (intrusive) and partially data-driven (non-intrusive). The numerical results allow us to assess the efficiency, and stability of the ROM, demonstrating its ability to accurately predict velocity and pressure fields, as well as quantities of interest, even when extrapolating in time.

Two main conclusions can be drawn from this study. First, when performing interpolation, i.e., remaining within the same parameter window for both physical and temporal variables, both the standard POD-ROM and the Hybrid POD-ROM achieve comparable levels of accuracy. In this case, the standard model incurs a slightly higher computational cost. Second, in the case of temporal extrapolation, both models again reach similar error levels; however, the Hybrid POD-ROM requires slightly more computational time. Overall, the main advantage of the standard POD-ROM lies in its lower computational cost during extrapolation, whereas the Hybrid POD-ROM offers greater flexibility for industrial applications. This is because it can be readily integrated with commercial software, as it does not rely on computing pressure boundary conditions at the full-order level.

Furthermore, the proposed POD-ROM framework could be successfully applied to more complex scenarios, such as three-dimensional flows, to enable faster simulations. This is particularly advantageous in many-query contexts, where numerical simulations must be performed for multiple values of physical parameters, for instance.
\\ \\
{\bf Acknowledgments:} This work has been supported by the Spanish Government Project PID2021-123153OB-C21. The research of M. Aza\"iez has been supported by  “Visiting Scholars” du “VII Plan Propio de Investigación y Transferencia” of the University of Seville (Spain) and the French National Research Agency ANR-22-CE46-0005 (NumOPT)\\

\bibliographystyle{abbrv}
\bibliography{reference}
\end{document}